\documentclass{amsart}
\usepackage{amssymb} 
\usepackage{color}
\usepackage[all]{xy}
\usepackage{enumerate}
\usepackage{ulem}

 \newtheorem{theorem}{Theorem}[section]

 \newtheorem{proposition}[theorem]{Proposition}
 \newtheorem{lemma}[theorem]{Lemma}
 \newtheorem{corollary}[theorem]{Corollary}
 \newtheorem{remark}[theorem]{Remark}
 \newtheorem{definition}[theorem]{Definition}
 \newtheorem{notation}[theorem]{Notation}
 \newtheorem{example}[theorem]{Example}

 \numberwithin{equation}{section}
 \newtheorem{warning}[theorem]{Warning}

 \def\subrel#1#2{\mathrel{\mathop{#2}\limits_{#1}}}

\def\cat{{\rm Cat}_{(\infty,1)}}

\def\map#1{{\rm Map}_{#1}}
\def\mapsc{{\rm Map}^{\rm sc}}

\def\twr{{\rm Tw}^r}
\def\btwr{\mathrm{TW}^r}
\def\bbtwr{\mathbb{T}\mathrm{W}^r}
\def\bstwr{\mathbb{T}\mathrm{w}^r}

\def\mapsc{{\rm Map}^{\rm sc}}

\def\setsc{{\rm Set}_{\Delta}^{\rm sc}}

\def\bsetdelta{{\rm bSet}_{\Delta}}
\def\setjoy{{\rm Set}_{\Delta}^{\rm Joy}}
\def\setplus{{\rm Set}_{\Delta}^{+}}
\def\catdeltaplus{{\rm Cat}_{\Delta}^+}
\def\setkan{{\rm Set}_{\Delta}^{\rm Kan}}
\def\setdelta{{\rm Set}_{\Delta}}

\def\twscbullet{T(\bullet)}

\def\twscn{T(n)}
\def\twscone{T(1)}

\def\twsccart{T_{\rm cart}}
\def\twscncart{T(n)_{\rm cart}}

\begin{document}

\title{A simple model for
twisted arrow $\infty$-categories}
\author{Takeshi Torii}
\address{Department of Mathematics, 
%Faculty of Sience, 
Okayama University,
Okayama 700--8530, Japan}
\email{torii@math.okayama-u.ac.jp}
%\thanks{}
%\thanks{The author was partially supported by
%JSPS KAKENHI Grant Numbers JP23K03113.}


\date{June 6, 2026 (version~2.2)}

\begin{abstract}
Twisted arrow $\infty$-categories
of $(\infty,1)$-categories
were introduced by Lurie,
and they have various applications
in higher category theory.
Abell\'{a}n Garc\'{i}a and Stern
gave a generalization to twisted arrow $\infty$-categories
of $(\infty,2)$-categories.
In this paper
we introduce another simple model
for twisted arrow $\infty$-categories
of $(\infty,2)$-categories.

\end{abstract}

\maketitle

\section{Introduction}

%{\color{red}
%There are many models of $(\infty,1)$-category.
%A quasi-category is one of the models of $(\infty,1)$-category.
%In this paper an $(\infty,1)$-category
%or an $\infty$-category means a quasi-category.
%We also use a complete Segal space
%which is also a model of $(\infty,1)$-category.
%}

The Yoneda embedding is a fundamental tool in category theory.
In higher category theory
the $(\infty,1)$-categorical
Yoneda embedding was proved
by Lurie \cite[Proposition~5.1.3.1]{Lurie1}.
There is a fully faithful functor
$C\to {\rm Fun}(C^{\rm op},\mathcal{S})$
from an $\infty$-category $C$
to the $\infty$-category of
presheaves on $C$
with values in the $\infty$-category $\mathcal{S}$
of spaces.
The $\infty$-category $\mathcal{S}$
is a full subcategory of the $\infty$-category
of (small) $(\infty,1)$-categories $\cat$.
In general,
it is difficult to work with functors into
$\mathcal{S}$ or $\cat$.
The $(\infty,1)$-categorical
Grothendieck construction
overcomes this difficulty.
Combining
\cite[Theorem~2.2.1.2]{Lurie1}
with
\cite[Proposition~4.2.4.4]{Lurie1},
it gives an equivalence between
the $\infty$-category
${\rm Fun}(C^{\rm op},\mathcal{S})$
of presheaves on $C$
and the $\infty$-category of right fibrations
over $C$.
%cf.~see also \cite[Corollary~A.32]{GHN}
It also gives an equivalence between
the $\infty$-category
${\rm Fun}(C^{\rm op},\cat)$
of functors
%from the opposite category
%$\mathcal{C}^{\rm op}$
into $\cat$ 
and the $\infty$-category of cartesian fibrations
over $C$
%(\cite[Theorem~2.2.1.2]{Lurie1}).
by combining \cite[Theorem~3.2.0.1]{Lurie1}
with \cite[Proposition~4.2.4.4]{Lurie1}.
It is easier to handle right fibrations
or cartesian fibrations
than functors into $\mathcal{S}$ or $\cat$.

%\bigskip

The twisted arrow $\infty$-category
of an $(\infty,1)$-category
is a right fibration associated
to the mapping space functor.
For an $\infty$-category $C$,
the twisted arrow $\infty$-category $\twr(C)$
was introduced by Lurie in \cite[\S5.2.1]{Lurie2}.
It comes equipped with a map $\twr(C)\to
C\times C^{\rm op}$
of $\infty$-categories
which is a right fibration classified by
the mapping space functor
${\rm Map}_{C}(-,-):
C^{\rm op}\times C\to \mathcal{S}$.

Twisted arrow $\infty$-categories
of $(\infty,1)$-categories play a fundamental role
in higher category theory.
They are used
to show
an existence of a dual equivalence
between two $(\infty,1)$-categories,
where a dual equivalence between $(\infty,1)$-categories
$C$ and $D$ means an equivalence
of $C$ and $D^{\rm op}$.
A pairing of $\infty$-categories between 
$C$ and $D$ 
is a functor
$C\times D\to \mathcal{S}$.
As mentioned above,
it is convenient to identify it with
a right fibration over $C^{\rm op}\times D^{\rm op}$.
As a matter of fact,
in \cite[Definition~5.2.1.5]{Lurie2},
Lurie has defined a pairing of $\infty$-categories 
to be a triple $(C,D,\lambda)$,
where $\lambda: M\to C\times D$
is a right fibration between $\infty$-categories.
A pairing is said to be perfect
if it is equivalent to a twisted arrow
$\infty$-category as right fibrations
under the equivalence
between ${\rm Fun}(C\times D,\mathcal{S})$
and the $\infty$-category of right fibrations
over $C^{\rm op}\times D^{\rm op}$.
If a pairing between $C$ and $D$ is perfect,
then it gives a dual equivalence
between $C$ and $D$.
This technique is
used to study iterated bar constructions
and Koszul duality for $\mathbb{E}_k$-algebras in
$\mathbb{E}_k$-monoidal $(\infty,1)$-categories in \cite[\S5.2]{Lurie2}.
%This is based on the fact that 
%the construction of twisted arrow $\infty$-catogories
%is symmetric monoidal.
It is also used in \cite{Torii2}
to construct a dual equivalence between
the $(\infty,1)$-category of monoidal
$(\infty,1)$-categories
and left adjoint oplax monoidal functors
and
%the $(\infty,1)$-category
that of monoidal $(\infty,1)$-categories
and right adjoint lax monoidal functors.

%\bigskip

Twisted arrow $\infty$-categories are also used
for other constructions
in the setting of higher category theory.
The $(\infty,1)$-categorical
(co)ends (\cite{Glasman}) and weighted (co)limits (\cite{GHN})
are constructed as some (co)limits
indexed by twisted arrow $\infty$-categories.
For an $\infty$-category $C$,
the $\infty$-category of spans in $C$
is constructed from some functors
from the twisted arrow $\infty$-categories of
the ordinary categories $[n]$
into $C$.
The construction of $\infty$-categories of spans
was introduced in \cite{Barwick2}.
It is shown that
the construction of twisted arrow $\infty$-categories
is a left adjoint to
the construction of $\infty$-categories of spans
in \cite[Theorem~A]{HHLN2}.

A generalization to the construction
of $\infty$-categories of twisted arrows
in $(\infty,2)$-categories
was introduced by Abell\'{a}n Garc\'{i}a and Stern
in \cite{Garcia-Stern}.
%For a fibrant scaled simplicial set $Z$,
For an $\infty$-bicategory $Z$,
%which is a model of an $(\infty,2)$-category,
they constructed an $\infty$-category $\twr(Z)$
equipped with a map
$p: \twr(Z)\to u_1Z\times u_1Z^{\rm op}$
of $\infty$-categories
which is a cartesian fibration classified by
the restricted mapping $\infty$-category functor
$\map{Z}(-,-):
(u_1Z\times u_1Z^{\rm op})^{\rm op}\simeq
u_1Z^{\rm op}\times u_1Z\to \cat$,
where $u_1Z$ is the underlying
$\infty$-category of $Z$.

The goal of this paper is 
to give another simple model for
twisted arrow $\infty$-categories of $(\infty,2)$-categories.
We first construct a complete Segal space $\btwr(Z)$
which is a lifting of the $\infty$-category $\twr(Z)$ 
with respect to
the Quillen equivalence
$i_1^*: \bsetdelta^{\rm CSS}\to \setjoy$ (\cite{Joyal-Tierney})
from the category of bisimplicial sets 
equipped with complete Segal space model structure
$\bsetdelta^{\rm CSS}$
(\cite{Rezk})
to the category of simplicial sets
equipped with Joyal model structure $\setjoy$
(\cite{Joyal-1, Lurie1}).
%(\cite{Joyal-1}, \cite[\S2.2.5]{Lurie1}).
Although
it is known that there is a lifting of $\twr(Z)$
to a complete Segal space
by \cite[Corollary~7.17]{HHLN2},
we construct $\btwr(Z)$
by making use of
combinatorial structure on scaled simplicial sets.

The construction of $\twr(Z)$ in \cite{Garcia-Stern}
uses a cosimplicial scaled simplicial set  
$Q(\bullet)$,
where the underlying simplicial set of
the $n$th scaled simplicial set $Q(n)$ is isomorphic to 
the $(2n+1)$-dimensional simplex $\Delta^{2n+1}$.
We introduce a cosimplicial scaled simplicial set $T(\bullet)$,
where
%the $n$th scaled simplicial set
$T(n)$ is a subcomplex of $Q(n)$
of dimension $(n+1)$.
We construct a bisimplicial set
$\bbtwr(Z)$
in the same way as $\btwr(Z)$
by using $T(\bullet)$ instead of $Q(\bullet)$.
We show that $\bbtwr(Z)$
is a complete Segal space
by comparing it with $\btwr(Z)$.
By applying the right Quillen equivalence $i_1^*$
to $\bbtwr(Z)$,
we obtain an $\infty$-category $\bstwr(Z)$.
As in $\twr(Z)$,
we show that 
the $\infty$-category $\bstwr(Z)$
comes equipped with a cartesian fibration
$p_{\mathbb{T}}: \bstwr(Z)\to u_1Z\times u_1Z^{\rm op}$.

The following is the main theorem of this paper,
which says that
the $\infty$-category $\bstwr(Z)$
equipped with the cartesian fibration $p_{\mathbb{T}}$
is equivalent to the twisted arrow $\infty$-category
$\twr(Z)$ with $p$.

\begin{theorem}
[{cf.~Theorem~\ref{thm:tw-equivalence-cartesian}}]  
\label{thm:main}
For any $\infty$-bicategory $Z$,
there is a natural equivalence of $\infty$-categories
\[ \twr(Z)\stackrel{\simeq}{\longrightarrow}
   \bstwr(Z)  \]
which makes the following diagram commute
\[ \xymatrix{
  \twr(Z)\ar[rr]\ar[dr]_p &&
  \bstwr(Z)\ar[dl]^{p_{\mathbb{T}}}\\
     &u_1Z\times u_1Z^{\rm op}. & }\]   
\end{theorem}

There are related constructions of twisted
arrow $\infty$-categories.
In \S\ref{section:comparison}
we discuss relationships between
them and the twisted arrow $\infty$-category in this paper.
In particular,
we compare twisted $2$-cell $\infty$-bicategories
of \cite{HNP} with the twisted arrow $\infty$-categories
of $(\infty,2)$-categories in this paper.
A twisted $2$-cell $\infty$-bicategory
is constructed from a fibrant marked simplicial category,
which is a model of an $(\infty,2)$-category,
and it comes equipped with
a scaled cocartesian fibration 
which is classified by the mapping $\infty$-category
functor of the marked simplicial category composed with
the twisted arrow construction of marked simplicial sets.
Although the construction of twisted arrow
$\infty$-category of an $(\infty,2)$-category 
in this paper
is less general than that of twisted $2$-cell
$\infty$-bicategory,
it is desirable to have a simple model of
twisted arrow $\infty$-category
%of an $(\infty,2)$-category
which is directly described
in terms of $\infty$-bicategories,
not through marked simplicial categories.

We will discuss an application of
the result of this paper in \cite{Torii5}.
We will introduce a notion of map monoidales
in monoidal $(\infty,2)$-categories,
and show that
the endomorphism $\infty$-category of a map monoidale
admits the structure of a duoidal $\infty$-category. 
We will use the main theorem in this paper
to show that one of the two monoidal structures
on the duoidal $\infty$-category is equivalent 
to convolution product.

%\bigskip

The organization of this paper is as follows:
In \S\ref{section:scaled-ss}
we review the theory of scaled simplicial sets.
We recall the notion of scaled anodyne maps
and the model structure on the category
of scaled simplicial sets.
We prove a lemma which gives a sufficient condition
for a map of scaled simplicial sets being scaled anodyne
for later use.
%In \S\ref{subsection:scaled-simplicial-sets}
%we recall the theory of scaled simplicial sets,
%which is a model of the theory of $(\infty,2)$-categories.
%We describe the model structure
%on the category of scaled simplicial sets
%introduced in \cite{Lurie3}.
%In \S\ref{subsection:scaled-anodyne-maps}
%we study scaled anodyne maps.
%We will prove a lemma
%which gives a sufficient condition
%on maps of scaled simplicial sets
%to be scaled anodyne. 
%
In \S\ref{section:twisted-arrow-categories}
we study the twisted arrow $\infty$-categories
of $(\infty,2)$-categories.
%In \S\ref{subsection:twisted-arrow-categories}
%we review the construction of
%twisted arrow $\infty$-category
%based on scaled simplicial sets
%introduced in \cite{Garcia-Stern}.
%The construction
%of twisted arrow $\infty$-categories
%determines a functor from the category of $\infty$-bicategories
%to the category of $\infty$-categories.
%Haugseng-Hebestreit-Linskens-Nuiten~\cite{HHLN2}
%showed that 
%this functor can be lifted to
%a functor with values in
%the category of complete Segal spaces.
%We will give another combinatorial proof of this fact.
%In \S\ref{subsection:lifting-Segal-spaces}
We construct a functor
from the category of fibrant scaled simplicial sets
to the category of complete Segal spaces,
which is a lifting of the construction
of twisted arrow $\infty$-categories in \cite{Garcia-Stern}.
%In \S\ref{subsection:twrZ-complete-Q}
%we show that it actually lands in 
%the category of complete Segal spaces.
%
In \S\ref{section:model-twisted-arrow}
we construct another model of the twisted
arrow $\infty$-category of an $(\infty,2)$-category.
First,
we introduce the cosimplicial scaled simplicial set
$\twscbullet$.
For a fibrant scaled simplicial set $Z$,
we define the bisimplicial set $\bbtwr(Z)$
by using $\twscbullet$,
and show that it is a complete Segal space.
Next,
we construct a cartesian fibration $p_{\mathbb{T}}$ and
finally prove the main theorem
(Theorem~\ref{thm:tw-equivalence-cartesian}).
In \S\ref{section:comparison}
we compare the twisted arrow $\infty$-category
in this paper with other related constructions.
In \S\ref{section:extension-TWr}
we construct a right Quillen functor
from the model category of scaled simplicial sets
to the complete Segal space model category,
which is equivalent to the functor
${\rm TW}^r(-)$ when restricted
to the full subcategory of $\infty$-bicategories.

%\bigskip

\begin{warning}\rm
Throughout this paper
we use the term $\infty$-category
%or $(\infty,1)$-category
for quasi-category.
%We regard ordinary $1$-categories as
%a special kind of $\infty$-categories,
%and we do not distinguish them with
%their nerves notationally.
\end{warning}

%\bigskip

\noindent
{\bf Acknowledgments}.
The author would like to thank
the referee for his/her careful reading,
valuable comments and suggestions.
The author was partially supported by
JSPS KAKENHI Grant Numbers JP23K03113.

%\newpage
%\input{notation}
\section{Notation}

%Let $\Delta$ be the category whose objects are finite,
%nonempty, totally orderd sets
%\[ [n]=\{0<1<\cdots <n\} \]
%for $n\ge 0$ and morphisms are order-preserving functions.
%Simplicial sets are functors from
%$\Delta^{\rm op}$ to the category of sets.
%Morphisms of simplicial sets are natural transformations.
For a nonempty finite totally ordered set $I$,
we denote by $\Delta^I$ the $(|I|-1)$-dimensional simplex
with $I$ as the set of vertices.
For $i\in I$,
$\Delta^{I-{\{i\}}}$ 
is the codimension $1$ face of $\Delta^I$
opposite to the vertex $i$.
For a subset $\emptyset \neq J\varsubsetneqq I$,
we set
\[ \Lambda^I_J=\bigcup_{i\in I-J}\Delta^{I-\{i\}}. \]
For $[r]=\{0<1<\cdots<r\}$,
we write $\Delta^r$ and $\Lambda^r_i$ for
$\Delta^{[r]}$ and $\Lambda^{[r]}_{\{i\}}$,
respectively,
as usual. 
We write $x\in X$
for a $0$-simplex $x$ of a simplicial set $X$.

We denote by $\setdelta$ 
the category of simplicial sets.
We can consider several model structures
on $\setdelta$ and its variants.
We denote by
$\setkan$ and $\setjoy$  
the category of simplicial sets equipped with
the Kan and Joyal model structures,
%(cf.~\cite[\S2.2.5]{Lurie1}),
respectively.
We denote by 
$\setplus$
the category of marked simplicial sets
equipped with the cartesian model
structure
(\cite[\S3.1.3]{Lurie1}).
We denote by $\setsc$
the category of scaled simplicial sets
equipped with the model structure
given by \cite[Theorem~4.2.7]{Lurie3}.
%We denteo by
%$\bsetdelta$
%the category of bisimplicial sets
%equipped with Reedy model structure.

We denote by $\bsetdelta$
the category of bisimplicial sets.
We write $\bsetdelta^{\rm Reedy}$
for the category of bisimplicial sets
equipped with the Reedy model structure.
We write $\bsetdelta^{\rm CSS}$
for the category of bisimplicial sets
equipped with the complete Segal space
model structure (\cite[Theorem~7.2]{Rezk}).

For a model category $M$,
we denote by
$M^{\circ}$
the full subcategory of $M$
spanned by those objects
which are both fibrant and cofibrant.
We write $M_{\infty}$ for
the underlying $\infty$-category of the model category $M$.

%%%\newpage
%%%\input{overview}
%%%\newpage
%\input{scaled-ss}
\section{Scaled simplicial sets and $\infty$-bicategories}
\label{section:scaled-ss}

In this section
we review scaled simplicial sets
and $\infty$-bicategories.
In \S\ref{subsection:scaled-simplicial-sets}
we recall the theory of scaled simplicial sets,
which is a model of the theory of $(\infty,2)$-categories.
We describe the model structure
on the category of scaled simplicial sets
introduced in \cite{Lurie3}.
In \S\ref{subsection:scaled-anodyne-maps}
we study scaled anodyne maps.
We prove a lemma
which gives a sufficient condition
on maps of scaled simplicial sets
to be scaled anodyne. 

\subsection{Scaled simplicial sets}
\label{subsection:scaled-simplicial-sets}

There are many models for $(\infty,2)$-categories.
A model of scaled simplicial sets
is one of them,
which was introduced by Lurie in \cite{Lurie3}.
In \cite{Barwick-SP} Barwick and Schommer-Pries 
introduced a theory of $(\infty,n)$-categories
by axiomatizing expected properties of
a model for $(\infty,n)$-categories,
and proved a uniqueness of theories
of $(\infty,n)$-categories up to equivalences.
Furthermore,
they verified that many known models satisfy the axioms
and, in particular, showed that  
almost all models for $(\infty,2)$-categories
are equivalent except for
a model of saturated $2$-trivial complicial sets. 
Finally, Gagna, Harpaz, and Lanari \cite{GHL}
showed that all known models for $(\infty,2)$-categories
are equivalent by constructing
a Quillen equivalence between the model category
of scaled simplicial sets
and that of saturated $2$-trivial complicial sets. 
Therefore,
the underlying homotopy theory
of the model category of scaled simplicial sets
is equivalent to
the theory of $(\infty,2)$-categories.
In this subsection
we recall the model structure
on the category of scaled simplicial sets,
which was constructed in \cite[Theorem~4.2.7]{Lurie3}.

First,
we recall the definition for scaled simplicial sets.
A scaled simplicial set $X$ is a pair
$(\overline{X},T_X)$
of simplicial set $\overline{X}$
and a set $T_{X}$ of $2$-dimensional simplices of $\overline{X}$
that contains all degenerate ones.
We say that a $2$-dimensional simplex of $\overline{X}$
is thin if it belongs to $T_X$.
A map of scaled simplicial sets
$f: X\to Y$
%$f: (\overline{X},T_X)\to (\overline{Y},T_Y)$
is a map of simplicial sets
$\overline{f}: \overline{X}\to \overline{Y}$
that satisfies $\overline{f}(T_X)\subset T_Y$.
We denote by
$\setsc$ the category of scaled simplicial sets.

For a simplicial set $S$,
we have two canonical scaled simplicial sets
$S_{\sharp}$ and $S_{\flat}$.
The scaled simplicial set $S_{\sharp}$
has the underlying simplicial set $S$ 
equipped with all $2$-dimensional
simplices as thin triangles.
On the other hand,
the scaled simplicial set $S_{\flat}$
has the underlying simplicial set $S$
equipped with degenerate $2$-dimensional
simplices as thin triangles.

Now,
we recall scaled anodyne maps of scaled simplicial sets
which characterize fibrant objects
in $\setsc$.
For a set $T$
of $2$-dimensional simplices of a simplicial set $\overline{X}$
and a subcomplex $\overline{Y}\subset \overline{X}$,
we write $T|_{\overline{Y}}$
for the subset $T\cap \overline{Y}_2$ of $T$.

\begin{definition}[{\cite[Definition~3.1.3]{Lurie3}}]
%[Definition of scaled anodyne map]
\rm
A collection of scaled anodyne maps is
the weakly saturated class  
of morphisms of scaled simplicial sets
generated by the following maps:

\begin{enumerate}

\item[An1:]
the inclusion maps
\[ (\Lambda^n_i,\{\Delta^{\{i-1,i,i+1\}}\}|_{\Lambda^n_i}\cup
   \{\mbox{\rm degenerate}\})
   \to (\Delta^n,\{\Delta^{\{i-1,i,i+1\}}\}\cup
   \{\mbox{\rm degenerate}\}) \]
for $n\ge 2$ and $0<i<n$,
  
\item[An2:]
the inclusion map
\[ (\Delta^4,T)\to
(\Delta^4,T\cup \{\Delta^{\{0,3,4\}},\Delta^{\{0,1,4\}}\}), \]
where
\[ T=\{\Delta^{\{0,2,4\}},\Delta^{\{1,2,3\}},
       \Delta^{\{0,1,3\}},\Delta^{\{1,3,4\}},
       \Delta^{\{0,1,2\}}\}\cup
       \{\mbox{\rm degenerate}\} , \]

\item[An3:]
the inclusion maps
\[ (\Lambda^n_0\coprod_{\Delta^{\{0,1\}}}\Delta^0,
    T|_{\Lambda^n_0\coprod_{\Delta^{\{0,1\}}}\Delta^0})
    \to
   (\Delta^n\coprod_{\Delta^{\{0,1\}}}\Delta^0, T)   \]
    for $n\ge 3$,
    where
\[ T= \{\Delta^{\{0,1,n\}}\}|_{\Delta^n\coprod_{\Delta^{\{0,1\}}}\Delta^0}
    \cup\{\mbox{\rm degenerate}\}.\]
\end{enumerate}
\end{definition}

%Next,
%we will define fibrant scaled simplicial sets.

\begin{definition}[{cf.~\cite[Theorem~4.2.7]{Lurie3}}]
%[Definition of fibrant scaled simplicial sets]
\label{def:definition-infinity-bicategory}
\rm
A scaled simplicial set is said to be
$\infty$-bicategory if
it admits extensions along every scaled anodyne map.
\end{definition}

\begin{remark}\rm
In \cite[Definition~4.1.1]{Lurie3}
scaled simplicial sets
in Definition~\ref{def:definition-infinity-bicategory}
are referred to as weak $\infty$-bicategories,
and the term $\infty$-bicategory was reserved
for fibrant objects in $\setsc$.
However,
it is shown that these two notions coincide
in \cite[Theorem~5.1]{GHL}.
\end{remark}

The model structure on the category
$\setsc$ of scaled simplicial sets
is given as follows:

\begin{definition}[{cf.~\cite[Theorem~4.2.7]{Lurie3},
\cite[Theorem~5.1]{GHL}}]
%[Definition of model structure  on scaled simplicial sets]
\rm
There exists a model structure on the category $\setsc$
of scaled simplicial sets whose
cofibrations are the monomorphisms and
whose fibrant objects are the $\infty$-bicategories.
We call a weak equivalence of scaled simplicial sets
a bicategorical equivalence.
This model structure is left proper and combinatorial.
\end{definition}

Now,
we study mapping spaces of $\setsc$.
For this purpose,
we recall 
relationships between
$\setsc$ and other model categories.
%and a mapping space on $\setsc$.
We will construct a diagram 
\[ \setkan 
   \stackrel{(-)^{\sharp}}{\subrel{u_0}{\rightleftarrows}}
   \setplus 
   \stackrel{(-)^{\flat}}{\subrel{U}{\leftrightarrows}}
   \setjoy 
   \stackrel{(-)_{\sharp}}{\subrel{u_1}{\rightleftarrows}}
   \setsc     \]
of Quillen adjunctions
between model categories,
where the upper arrows are left adjoint
and the lower arrows are right adjoint.

First,
we consider the pair $((-)^{\sharp}, u_0)$ of functors.
The functor $(-)^{\sharp}: \setkan\to\setplus$
assigns to a simplicial set $X$
the marked simplicial set $X^{\sharp}=(X,X_1)$.
Its right adjoint $u_0: \setplus\to \setkan$
assigns to a marked simplicial set $(X,S)$
the subcomplex of $X$ spanned by those simplices
whose edges are all marked.
In particular,
when $Z$ is a fibrant marked simplicial set,
$u_0Z$ is a Kan complex.
We call $u_0Z$ the underlying $\infty$-groupoid of $Z$.

We will show that the pair $((-)^{\sharp}, u_0)$ is a Quillen adjunction
(cf.~\cite[Example~A.3.2.22]{Lurie1}).
It is clear that $(-)^{\sharp}$ preserves
cofibrations.
Let $i: A\to B$ be a trivial cofibration
in $\setkan$.
We shall show that $i^{\sharp}: A^{\sharp}\to B^{\sharp}$
is a trivial cofibration in $\setplus$.
For an $\infty$-category $Y$,
we denote by $Y^{\natural}$
the fibrant marked simplicial set
whose underlying simplicial set $Y$
equipped with equivalences as marked edges.
We recall that $\setplus$ is a simplicial model category
with mapping object $\map{}^{\sharp}(-,-)$
by \cite[Corollary~3.1.4.4]{Lurie1}.
By the definition of weak equivalence in $\setplus$,
it suffices to show that the induced map
on mapping spaces
$\map{}^{\sharp}(B^{\sharp},Y^{\natural})\to
\map{}^{\sharp}(A^{\sharp},Y^{\natural})$
is a trivial Kan fibration
for any $\infty$-category $Y$.
%for any fibrant marked simplicial set $Z=(\overline{Z},E_Z)$.
This follows by observing that
$\map{}^{\sharp}(C^{\sharp},Y^{\natural})$ is isomorphic to 
$\map{}(C,Y^{\simeq})$
for a simplicial set $C$,
where $Y^{\simeq}$
is the underlying Kan complex
%$\infty$-groupoid
of $Y$.

Next,
we consider the pair $((-)^{\flat},U)$ of functors.
The functor $(-)^{\flat}:\setjoy\to \setplus$
assigns to a simplicial set $X$
the marked simplicial set
$X^{\flat}=(X,s_0(X_0))$.
%$X^{\flat}=(X,\{\mbox{\rm degenerate}\})$.
Its right adjoint $U$ is the forgetful functor.
The pair $((-)^{\flat},U)$
is a Quillen equivalence
by \cite[Proposition~3.1.5.3]{Lurie1}.

Finally,
we consider the pair $((-)_{\sharp},u_1)$
of functors.
The functor $(-)_{\sharp}: \setjoy\to \setsc$
is given by $X\mapsto X_{\sharp}$.
Its right adjoint $u_1: \setsc\to \setjoy$
assigns to a scaled simplicial set $Y=(\overline{Y},T_Y)$
the subcomplex of $\overline{Y}$
spanned by those simplices whose
$2$-dimensional faces are all thin.
By \cite[Remark~4.1.3]{Lurie3},
$u_1Z$ is an $\infty$-category
when $Z$ is an $\infty$-bicategory.
We call $u_1Z$ the underlying $\infty$-category
of $Z$. 

The model structure on 
$\setsc$ is cartesian closed
by \cite[Proposition~3.1.8 and Lemma~4.2.6]{Lurie3}
(see, also,
\cite[Remark~1.31]{GHL} or
\cite[the paragraph before Lemma~1.22]{GHL2}).
Thus,
we have a function object ${\rm FUN}(A,B)$
in $\setsc$
for scaled simplicial sets $A$ and $B$.
When $Z$ is an $\infty$-bicategory,
${\rm FUN}(A,Z)$ is also an $\infty$-bicategory,
and we denote
by ${\rm Fun}(A,Z)$
its underlying $\infty$-category
$u_1{\rm FUN}(A,Z)$.
In particular,
${\rm Fun}(A,Z)^{\natural}$
is a fibrant marked simplicial set
for an $\infty$-bicategory $Z$.
By applying the functor $u_0$,
we obtain a Kan complex
\[ \mapsc(A,Z)=u_0{\rm Fun}(A,Z)^{\natural}. \]
Note that
$\mapsc(A,Z)$ is the underlying $\infty$-groupoid
of the $\infty$-category ${\rm Fun}(A,Z)$
since $u_0Y^{\natural}=Y^{\simeq}$
for any $\infty$-category $Y$.
%where $Y^{\simeq}$ is the underlying $\infty$-groupoid
%of $Y$.

We will prove that the pair $((-)_{\sharp},u_1)$
is a Quillen adjunction.
It is clear that $(-)_{\sharp}$ preserves cofibrations.
Thus,
it suffices to show that
$(-)_{\sharp}$ preserves trivial cofibrations.
Let $i: A\to B$ be a trivial cofibration
in $\setjoy$.
We shall show that $i_{\sharp}: A_{\sharp}\to B_{\sharp}$
is a trivial cofibration in $\setsc$.
By \cite[Lemma~1.22]{GHL2},
it suffices to show that the induced map
on mapping spaces
$\mapsc(B_{\sharp},Z)\to\mapsc(A_{\sharp},Z)$
is a trivial Kan fibration
for any $\infty$-bicategory $Z$.
This follows by observing that
$\mapsc(C_{\sharp},Z)$ is isomorphic to 
$\map{}^{\sharp}(C^{\flat},(u_1Z)^{\natural})$
for any simplicial set $C$.

We need the following proposition
in \S\ref{subsection:lifting-Segal-spaces}
below.

%\begin{corollary}\label{cor:cof-induce-map-sc-fibration}
\begin{proposition}\label{prop:cof-induce-map-sc-fibration}
Let $X$ be an $\infty$-bicategory
and let $A\to B$ be a cofibration in $\setsc$.
The induced map
$\mapsc(B,X)\to \mapsc(A,X)$
is a Kan fibration of Kan complexes
in $\setkan$.
\end{proposition}

In order to prove Proposition~\ref{prop:cof-induce-map-sc-fibration},
we need a preliminary lemma.
The construction $Y\mapsto Y^{\natural}$
determines a functor
\[ (-)^{\natural}: (\setjoy)^{\circ}
   \longrightarrow (\setplus)^{\circ}. \]
By \cite[Proposition~3.1.3.5]{Lurie1},
the functor $(-)^{\natural}$ preserves and
reflects weak equivalences.
We show that
$(-)^{\natural}$ also preserves fibrations.

\begin{lemma}
\label{lemma:categorical-fibration-implies-Kan-fibration}
For a categorical fibration $p: X\to Y$ 
of $\infty$-categories,
the induced map
$p^{\natural}: X^{\natural}\to Y^{\natural}$
is a fibration of fibrant objects in
$\setplus$.
\end{lemma}

\proof
Let $i: A\to B$
be a trivial cofibration in $\setplus$.
We consider a commutative diagram of solid arrows
\begin{equation}\label{eq:categorical-fib-to-fib-cartesian}
  \xymatrix{
     A\ar[r]^f\ar[d]_{i} & X^{\natural}\ar[d]^{p^{\natural}}\\
     B\ar[r]_g\ar@{..>}[ur]^{h}& Y^{\natural}.}
   \end{equation}
We need to show that there is a dotted arrow
$h: B\to X^{\natural}$ which makes
the whole diagram commute.

%We recall that $\map{}^{\sharp}(-,-)$
%is the mapping object of the simplicial model category
%$\setplus$ (cf.~\cite[Corollary~3.1.4.4]{Lurie1}).
We have a commutative diagram
\[ \xymatrix{
     {\rm Map}^{\sharp}(B,X^{\natural})
     \ar[r]^{p^{\natural}_*}\ar[d]_{i^*}&
     {\rm Map}^{\sharp}(B,Y^{\natural})\ar[d]^{i^*}\\
     {\rm Map}^{\sharp}(A,X^{\natural})
     \ar[r]^{p^{\natural}_*}&
     {\rm Map}^{\sharp}(A,Y^{\natural})\\
 }\]
of Kan complexes.
The vertical arrows are trivial Kan fibrations
by \cite[Lemma~3.1.3.6]{Lurie1}
and the definition of weak equivalences
in $\setplus$
(cf.~\cite[Proposition~3.1.3.3]{Lurie1}).
Thus,
we obtain
  a homotopy equivalence of Kan complexes
\begin{equation}\label{eq:weak-eq-Map-sharp}
  (i^*,p^{\natural}_*):
  {\rm Map}^{\sharp}(B,X^{\natural})
   \stackrel{\simeq}{\longrightarrow}
   {\rm Map}^{\sharp}(A,X^{\natural})
   \times_{{\rm Map}^{\sharp}(A,Y^{\natural})}           
         {\rm Map}^{\sharp}(B,Y^{\natural}).
         \end{equation}

Commutative
diagram~(\ref{eq:categorical-fib-to-fib-cartesian})
determines a vertex $(f,g)$
on the right hand side of (\ref{eq:weak-eq-Map-sharp}).
By homotopy equivalence~(\ref{eq:weak-eq-Map-sharp}),
we can take $k\in {\rm Map}^{\sharp}(B,X^{\natural})$
and a homotopy
$K=(K_1,K_2): \Delta^1\to
   {\rm Map}^{\sharp}(A,X^{\natural})
   \times_{{\rm Map}^{\sharp}(A,Y^{\natural})}           
         {\rm Map}^{\sharp}(B,Y^{\natural})$
from $(f,g)$ to $(i^*,p^{\natural}_*)(k)$.
By using $k$ and $K$,
we obtain a commutative diagram of solid arrows
\begin{equation}\label{eq:lifting-H}
  \xymatrix{
  M\ar[rr]^{\overline{K}_1}
  \ar@{^{(}->}[d] && X^{\natural}\ar[d]^{p^{\natural}}\\
  B\times(\Delta^{1})^{\sharp}
  \ar[rr]_{K_2}\ar@{..>}[urr]^H
  && Y^{\natural},}
  \end{equation}
where $M$ is a subobject
of $B\times (\Delta^1)^{\sharp}$
given by $(A\times(\Delta^1)^{\sharp})\coprod_{A\times \Delta^{\{1\}}}
(B\times \Delta^{\{1\}})$,
and
$\overline{K}_1$ is a morphism whose restrictions
are given by
$\overline{K}_1|_{A\times(\Delta^1)^{\sharp}}=K_1$ and
$\overline{K}_1|_{B\times \Delta^{\{1\}}}=k$.

Combining \cite[Proposition~3.1.2.3]{Lurie1}
with the fact that
$\Delta^{\{1\}}\to (\Delta^1)^{\sharp}$
is a marked anodyne morphism
by \cite[Definition~3.1.1]{Lurie1},
the inclusion map $M\to
B\times (\Delta^1)^{\sharp}$ is marked anodyne.
We can verify that $p^{\natural}$ 
satisfies the conditions in
\cite[Proposition~3.1.1.6]{Lurie1}
by using the fact that
$p: X\to Y$ is a categorical fibration
of $\infty$-categories,
%\cite[Proposition~1.2.4.3]{Lurie1},
the dual of \cite[Corollary~2.4.6.5]{Lurie1},
%we can extend $\overline{K}_1$
%on $M\cup (B^{(0)}\times (\Delta^1)^{\sharp})$,
%where $B^{(0)}$ is the $0$-skelton of $B$.
%Furthermore,
and \cite[Proposition~2.4.1.5]{Lurie1}.
%and \cite[Remark~2.4.1.4]{Lurie1}.
Hence,
$p^{\natural}$ 
has the right lifting property
with respect to all marked anodyne morphisms,
and we can construct a lifting
$H: B\times (\Delta^1)^{\sharp}\to X^{\natural}$
which makes whole diagram~(\ref{eq:lifting-H}) commute.
Then $h=H|_{B\times\Delta^{\{0\}}}$
gives the desired lifting of
commutative diagram~(\ref{eq:categorical-fib-to-fib-cartesian}).
\qed

%\bigskip

\proof[Proof of Proposition~\ref{prop:cof-induce-map-sc-fibration}]
The induced map ${\rm FUN}(B,X)\to {\rm FUN}(A,X)$
is a fibration of
fibrant objects in $\setsc$
since the model category
$\setsc$ is cartesian closed.
The fact that $u_1$ is a right Quillen functor
implies that
${\rm Fun}(B,X)\to {\rm Fun}(A,X)$
is a fibration of fibrant objects
in $\setjoy$.
The lemma follows
from Lemma~\ref{lemma:categorical-fibration-implies-Kan-fibration}.
\qed

\subsection{A lemma on scaled anodyne maps}
\label{subsection:scaled-anodyne-maps}

In this subsection
we prove a lemma which gives
a sufficient condition for an inclusion map
of scaled simplicial sets into a simplex 
being scaled anodyne.

First,
we recall the notation of simplicial sets.
For a nonempty finite totally ordered set $I$,
we denote by
$\Delta^I$ the $(|I|-1)$-dimensional simplex
with $I$ as the set of vertices.
For $i\in I$,
$\Delta^{I-{\{i\}}}$ 
is the codimension $1$ face of $\Delta^I$
opposite to the vertex $i$.
For a subset $\emptyset \neq J\varsubsetneqq I$,
we set
\[ \Lambda^I_J=\bigcup_{i\in I-J}\Delta^{I-\{i\}}. \]
For $[n]=\{0<1<\cdots<n\}$,
we write $\Delta^n$ and $\Lambda^n_i$ for
$\Delta^{[n]}$ and $\Lambda^{[n]}_{\{i\}}$,
respectively,
as usual. 

Let $\Delta^n_{\dag}=(\Delta^n,T_{\Delta^n_{\dag}})$
be a scaled simplicial set
whose underlying simplicial set is $\Delta^n$.
For a subcomplex $L$ of $\Delta^n$,
we denote by $L_{\dag}$
the scaled simplicial set
whose underlying simplicial set is $L$
equipped with the induced scaling.

We extensively use the following lemma. 

\begin{lemma}[{cf.~\cite[Lemma~1.18]{Garcia-Stern}}]
\label{lemma:simplex-scaled-anodyne}
Let $n\ge 3$
and let $M$ be a nonempty subset of
%$\{0,1,\ldots,n-1\}$.
 $[n]-\{n\}$.
We assume that
there exists an integer $s$ with
$0\le s<t$
such that
$s\not\in M$ and $a\in M$ for all $s<a\le t$,
where $t$ is the largest number of $M$.
Furthermore,
we assume
that $|M|\le n-2$ and
that the triangle $\Delta^{[n]-M}$
is not thin in $\Delta^n_{\dag}$
when $|M|=n-2$.
If triangles $\Delta^{\{b,t,t+1\}}$
are thin in $\Delta^n_{\dag}$ for all $s\le b< t$,
then the inclusion map
$\Lambda^n_{M,\dag}\to \Delta^n_{\dag}$
is scaled anodyne.
\end{lemma}

\proof
We notice that
all thin triangles in $\Delta^n_{\dag}$
are contained in $\Lambda^n_{M,\dag}$
by the assumption
that $|M|\le n-2$ and
$\Delta^{[n]-M}$ is not thin
when $|M|=n-2$.
We prove the lemma by induction
on the cardinality $|M|$ of $M$.

First,
we consider the case $|M|=1$.
We write $M=\{m\}$
with $0<m<n$.
In this case
we have $0\le s=m-1$ and $t=m<n$.
By the assumption,
$\Delta^{\{m-1,m,m+1\}}$ is thin in $\Delta^n_{\dag}$.
There is a pushout diagram
\[ \begin{array}{ccc}
  (\Lambda^n_m,\{\Delta^{\{m-1,m,m+1\}}\}
  \cup\{{\rm degenerate}\})
  & \longrightarrow &
  (\Delta^n, \{\Delta^{\{m-1,m,m+1\}}\}
  \cup\{{\rm degenerate}\})\\
  \bigg\downarrow & & \bigg\downarrow \\
  \Lambda^n_{m,\dag} & \longrightarrow &
  \Delta^n_{\dag}\\
   \end{array}\]
of scaled simplicial sets.
Since the top horizontal arrow 
is a scaled anodyne map, 
so is the bottom horizontal arrow.  

Next,
we suppose that $|M|>1$.
We let $m=\min M$ and set $M'=M-\{m\}$.
When we regard $M'$ as a subset of $[n]-\{n\}$,
we can verify that $M'$
satisfies the assumptions on the lemma
if $M$ does.
By the hypothesis of induction,
$\Lambda^{[n]}_{M',\dagger}\to \Delta^{[n]}_{\dag}$
is scaled anodyne.
Thus, it suffices to show that
the inclusion map
$\Lambda^{[n]}_{M,\dagger}\to
\Lambda^{[n]}_{M',\dagger}$
is scaled anodyne.
There is a pushout diagram
\[ \begin{array}{ccc}
  \Lambda^{[n]-\{m\}}_{M',\dag}&
  \longrightarrow  &
  \Delta^{[n]-\{m\}}_{\dag}\\
  \bigg\downarrow & & \bigg\downarrow \\[3mm]
  \Lambda^{[n]}_{M,\dag}&\longrightarrow &
  \Lambda^{[n]}_{M',\dag}\\  
   \end{array}\]
of scaled simplicial sets.
We have an isomorphism
$\theta: [n]-\{m\}\cong [n-1]$
of finite ordered sets.
We can verify that
the subset $\theta(M')\subset [n-1]-\{n-1\}$
satisfies the assumptions on the lemma
if $M\subset [n]-\{n\}$ does.
Thus,
the top horizontal arrow is scaled anodyne
by the hypothesis of induction,
and so is the bottom horizontal arrow.
\qed

\begin{remark}\rm
The dual statement of Lemma~\ref{lemma:simplex-scaled-anodyne}
also  holds:
Let $n\ge 3$
and let $M$ be a nonempty subset of
 $[n]-\{0\}$.
%$\{1,\ldots,n\}$.
We assume that
there exists an integer $s$ with
$t< s\le n$
such that
$s\not\in M$ and $a\in M$ for all $t\le a< s$,
where $t$ is the least number of $M$.
Furthermore,
we assume
that $|M|\le n-2$ and
that the triangle $\Delta^{[n]-M}$
is not thin in $\Delta^n_{\dag}$
when $|M|=n-2$.
If triangles $\Delta^{\{t-1,t,b\}}$
are thin in $\Delta^n_{\dag}$ for all $t<b\le s$,
then the inclusion map
$\Lambda^n_{M,\dag}\to \Delta^n_{\dag}$
is scaled anodyne.
\end{remark}

\begin{remark}\rm
When $|M|<n-2$,
we notice that Lemma~\ref{lemma:simplex-scaled-anodyne}
is a special case of
\cite[Lemma~1.18]{Garcia-Stern}
by setting $\mathcal{A}=\{\{a\}|\ a\in [n]-M\}$
with pivot point $t=\max{M}$.
In this case ${\rm Bas}(\mathcal{A})=\{Z\}$
with $Z=[n]-M$,
and
$l^Z_{t-1}=s$ and $l^Z_{t}=t+1$.
\end{remark}

\begin{example}\rm
For example,
we consider the case in which $n=5$ and $M=\{1,3,4\}
\subset [5]-\{5\}=\{0,1,2,3,4\}$.
In this case we have $s=2$ and $t=4$.
Suppose that $\Delta^5_{\dagger}$
is a scaled simplicial set such that
the $2$-simplex $\Delta^{\{0,2,5\}}$ is not thin.
If $\Delta^{\{2,4,5\}}, \Delta^{\{3,4,5\}}$ are thin
in $\Delta^5_{\dagger}$,
then the inclusion map $\Lambda^5_{M,\dagger}\to \Delta^5_{\dagger}$
is scaled anodyne.
\end{example}

%%%\newpage
%\input{GHmodel}
\section{Twisted arrow $\infty$-categories}
\label{section:twisted-arrow-categories}

In this section
we study the twisted arrow $\infty$-categories
for $\infty$-bicategories.
In \S\ref{subsection:twisted-arrow-categories}
we review the construction of
twisted arrow $\infty$-category
%based on scaled simplicial sets
introduced in \cite{Garcia-Stern}.
The construction
of twisted arrow $\infty$-categories
determines a functor from the category of $\infty$-bicategories
to the category of $\infty$-categories.
Haugseng-Hebestreit-Linskens-Nuiten~\cite[Corollary~7.17]{HHLN2}
showed that 
this functor can be lifted to
a functor with values in
the category of complete Segal spaces.
We give another combinatorial proof of this fact.
In \S\ref{subsection:lifting-Segal-spaces}
we construct a functor to the category of Segal spaces.
In \S\ref{subsection:twrZ-complete-Q}
we show that it lands in 
the category of complete Segal spaces.

\subsection{Twisted arrow $\infty$-categories}
\label{subsection:twisted-arrow-categories}

In this subsection
we recall the construction
of twisted arrow $\infty$-category
introduced in \cite{Garcia-Stern}.
For $n\ge 0$,
we consider the ordinary category
$[n]\star[n]^{\rm op}$,
where $(-)\star(-)$ is the join of 
ordinary categories (cf.~\cite[\S1.2.8]{Lurie1}).
By applying the nerve functor 
to $[n]\star[n]^{\rm op}$,
we obtain the simplicial set
$\Delta^n\star\Delta^{n,{\rm op}}$,
which is isomorphic to $\Delta^{2n+1}$.
We represent it as the following diagram
\[ \xymatrix{
  00\ar[r]\ar[d]
  & 01\ar[r]\ar[d]
  & \cdots\ar[r]
  & 0n\ar[d] \\
  10
  & 11\ar[l]
  & \cdots\ar[l]
  & 1n.\ar[l] \\
}   \]

\begin{definition}\rm
We recall
the scaled simplicial set $Q(n)
=(\Delta^n\star\Delta^{n,{\rm op}}, T_{Q(n)})$
in \cite[Definition~2.2]{Garcia-Stern}.
The underlying simplicial set
of $Q(n)$ is $\Delta^n\star\Delta^{n,{\rm op}}$
and the set of thin triangles $T_{Q(n)}$
is given by
\[ \begin{array}{rcll}
  T_{Q(n)}&=& &
  \{\Delta^{\{ik,ik',ik''\}}|\ i=0,1,\ 0\le k<k'<k''\le n\} \\
  & & \cup &
  \{\Delta^{\{0k,0k',1k''\}}|\ 0\le k< k'\le k''\le n\}\\
  & & \cup &
  \{\Delta^{\{1k,1k',0k''\}}|\ 0\le k< k'\le k''\le n\}
  \\  
            & & \cup &
  \{\mbox{\rm degenerate}\}.\\ 
\end{array}\]
\end{definition}

%\begin{definition}\rm
\begin{notation}\rm
For a subcomplex $K$ of $\Delta^n\star\Delta^{n,{\rm op}}$,
we denote by
% $K_{\ddagger}$
  $K_{\S}$
the scaled simplicial set whose underlying simplicial set
is $K$ equipped with the induced scaling from $Q(n)$.
\end{notation}
%\end{definition}

We notice that
the collection $Q(\bullet)=\{Q(n)\}_{n\ge n}$
forms a cosimplicial object
of scaled simplicial sets.
Furthermore,
the inclusion maps
$\Delta^n_{\sharp}\to Q(n)$ and
$\Delta^{n,{\rm op}}_{\sharp}\to Q(n)$
induce maps of
cosimplicial objects
$\Delta^{\bullet}_{\sharp}=\{\Delta^n_{\sharp}\}_{n\ge 0}\to
Q(\bullet)$
%\{Q(n)\}_{n\ge 0}$
and 
$\Delta^{\bullet,{\rm op}}_{\sharp}=
\{\Delta^{n,{\rm op}}_{\sharp}\}_{n\ge 0}\to
Q(\bullet)$
%\{Q(n)\}_{n\ge 0}$
in scaled simplicial sets.

%Let $Z=(\overline{Z},T_Z)$ be an $\infty$-bicategory.
%We regard $Z$ as a representative 
%of an $(\infty,2)$-category $\mathcal{Z}$.
%We denote by $u_1Z$
%the subcomplex of $\overline{Z}$
%spanned by those simplices whose
%$2$-dimensional faces are all thin.
%Then $u_1Z$ is an $\infty$-category and
%represents the underlying $(\infty,1)$-category
%$u_1\mathcal{Z}$ of $\mathcal{Z}$.

Now,
we recall the construction of
twisted arrow $\infty$-categories.

\begin{definition}[{\cite[Definition~2.4]{Garcia-Stern}}]
\label{def:quasi-category-twr-Qn}
\rm
Let $Z=(\overline{Z},T_Z)$ be an $\infty$-bicategory.
For $n\ge 0$,
we consider the set
${\rm Hom}_{\setsc}(Q(n),Z)$
of morphisms of scaled simplicial sets
from $Q(n)$ to $Z$.
This determines a simplicial set
\[ \twr(Z)=\{{\rm Hom}_{\setsc}(Q(n),Z)\}_{n\ge 0}. \]
The morphisms
$\Delta^{\bullet}_{\sharp}\to Q(\bullet)$
%$\{\Delta^n_{\sharp}\}_{n\ge 0}\to \{Q(n)\}_{n\ge 0}$
and
$\Delta^{\bullet,{\rm op}}_{\sharp}\to Q(\bullet)$
%$\{\Delta^{n,{\rm op}}_{\sharp}\}_{n\ge 0}
%\to \{Q(n)\}_{n\ge 0}$
induce a morphism of simplicial sets
\[ p: \twr(Z)\to u_1Z\times u_1Z^{\rm op}. \]
%\begin{theorem}[{\cite[Theorem~0.1]{Garcia-Stern}}]
which is a cartesian fibration 
by \cite[Theorem~0.1]{Garcia-Stern}.
In particular,
$\twr(Z)$ is an $\infty$-category
since $u_1Z\times u_1Z^{\rm op}$
is an $\infty$-category.
The cartesian fibration $p$ is classified
by the restricted mapping category functor
\[ \map{Z}(-,-):
   u_1Z^{\rm op}\times u_1Z
   \longrightarrow   \cat \]
by \cite[Theorem~0.1]{Garcia-Stern}.
%\[   \mathcal{Z}^{\mbox{\scriptsize\rm $1$-op}}\times
%   \mathcal{Z}\longrightarrow
%   \cat. \]
%where $\mathcal{Z}$ is the $(\infty,2)$-category
%represented by $Z$.
%for any $\infty$-bicategory $Z$.
%\end{theorem}
\end{definition}

\subsection{Liftings to Segal spaces
for twisted arrow $\infty$-categories}
\label{subsection:lifting-Segal-spaces}

The construction
$Z\mapsto \twr(Z)$
determines a functor
$\twr(-):(\setsc)^{\circ}\to (\setjoy)^{\circ}$.
Haugseng-Hebestreit-Linskens-Nuiten~\cite[Corollary~7.17]{HHLN2}
showed that 
this functor can be lifted to
a functor with values in
the category of complete Segal spaces.
We give another combinatorial proof of this fact.
In this subsection
we show that
the functor $\twr(-)$
lifts to the category of Segal spaces.
In the next subsection  
we show that
it takes values in the category
of complete Segal spaces.

\begin{definition}
\rm
Let $Z$ be an $\infty$-bicategory.
For $n\ge 0$,
we set
\[ \btwr(Z)_n=
   \mapsc(Q(n),Z)\in\setdelta.\]
The collection $\{\btwr(Z)_n\}_{n\ge 0}$
forms a bisimplicial set
$\btwr(Z)$.
This construction determines a functor
\[ \btwr(-):
  (\setsc)^{\circ}
  \longrightarrow \bsetdelta,\]
where $\bsetdelta$ is the category of bisimplicial sets.
\end{definition}

The goal of this subsection
is to prove the following proposition.

\begin{proposition}\label{prop:lifting-SS}
The bisimplicial set 
$\btwr(Z)$
is a Segal space
for any $\infty$-bicategory $Z$.
\end{proposition}

To prove Proposition~\ref{prop:lifting-SS},
we give some preliminary lemmas.
For the model category
$\setsc$ of scaled simplicial sets,
we can consider the Reedy model structure
on the category of cosimplicial objects
of $\setsc$
(see, for example,
\cite{Hirschhorn,Hovey}
for the Reedy model structure).
Let $A^{\bullet}$ be a cosimplicial object
of scaled simplicial sets.
We assume that $A^{\bullet}$
is cofibrant in the Reedy model structure
on the category of cosimplicial objects
of $\setsc$.
In other words,
the map $L_nA^{\bullet}\to A^n$ is a monomorphism
for all $n\ge 0$,
where $L_nA^{\bullet}$ is the $n$th latching object
of $A^{\bullet}$.
For $n\ge 0$,
we recall that the co-Segal map  
\begin{equation}\label{eq:co-Segal-map}
  \overbrace{A^1\subrel{A^0}{\coprod}\cdots
    \subrel{A^0}{\coprod}A^1}^n
  \simeq A^{\{0,1\}}\subrel{A^{\{1\}}}{\coprod}\cdots
  \subrel{A^{\{n-1\}}}{\coprod}A^{\{n-1,n\}}
  \to A^{\{0,1,\ldots,n\}}\simeq 
  A^n
\end{equation}
is induced by inert morphisms $[i]\to [n]$
for $i=0,1$.
We say that $A^{\bullet}$ satisfies the co-Segal condition
if map (\ref{eq:co-Segal-map})
is a weak equivalence for all $n\ge 0$.

\begin{lemma}\label{lemma:general-reedy-Segal-space}
Let $X$ be an $\infty$-bicategory.
If $A^{\bullet}$ is a Reedy cofibrant
cosimplicial object of $\setsc$,
then 
$\mapsc(A^{\bullet},X)$
is a Reedy fibrant simplicial object
of $\setkan$.
In addition,
if $A^{\bullet}$ satisfies the co-Segal condition,
then
$\mapsc(A^{\bullet},X)$
is a Segal space.
\end{lemma}

\proof
The first part follows 
by Proposition~\ref{prop:cof-induce-map-sc-fibration}.
If $A^{\bullet}$ satisfies the co-Segal condition,
then we see that $\mapsc(A^{\bullet},X)$
satisfies the Segal condition.
\qed

\bigskip

By using Lemma~\ref{lemma:general-reedy-Segal-space},
we will show that
the bisimplicial set $\btwr(Z)$ is Reedy fibrant.
%for any $\infty$-bicategory $Z$.
%For notational ease,
%we introduce a functor $Q: \setdelta\to \setsc$.

\begin{lemma}
\label{lemma:twr-Z-Segal-space}
The cosimplicial object
$Q(\bullet)$ of scaled simplicial sets
is Reedy cofibrant.
Hence,
the bisimplicial set
${\rm TW}^r(Z)$
is Reedy fibrant
for any $\infty$-bicategory $Z$.
%The bisimplicial set
%$\twr(Z)_{\bullet}$
%is Reedy fibrant
%for any $\infty$-bicategory $Z$.
\end{lemma}

\proof
The first part follows by observing that
the $n$th latching object of $Q(\bullet)$
is isomorphic to the subcomplex
$\cup_{i=0}^n(\Delta^{[n]-\{i\}}\star\Delta^{[n]-\{i\},{\rm op}})_{\S}$
of $Q(n)$.
The second part follows 
by Lemma~\ref{lemma:general-reedy-Segal-space}
since ${\rm TW^r(Z)}={\rm Map}^{\rm sc}(Q(\bullet),Z)$.
\qed

\bigskip

By \cite[Proposition~3.1.5]{Hovey},
the cosimplicial object $Q(\bullet)$ induces an adjunction
\[ L_Q: \setdelta\rightleftarrows \setsc : R_Q, \]
where $L_Q$ is 
the left Kan extension of $Q(\bullet): \Delta\to \setsc$
along the Yoneda embedding $\Delta\to \setdelta$.
%the coend
%$Q(\bullet)\otimes K$ for $K\in \setdelta$.
For notational ease,
we write $Q(K)$
for the scaled simplicial set $L_Q(K)$.
By \cite[Proposition~5.4.1]{Hovey}
and Lemma~\ref{lemma:twr-Z-Segal-space},
the left adjoint $L_Q: \setdelta\to\setsc$
preserves cofibrations.
%\begin{remark}\rm
Therefore,
for a monomorphism $A\to B$ of simplicial sets,
$Q(A)\to Q(B)$ is a monomorphism of scaled simplicial sets.  
In particular,
for a subcomplex $K$ of $\Delta^n$,
the scaled simplicial set $Q(K)$
is isomorphic to a subcomplex of $Q(n)$ given by
$\cup_{\Delta^I\subset K}
(\Delta^I\star \Delta^{I, {\rm op}})_{\S}$,
where the union ranges over subsets $I\subset [n]$
such that $\Delta^I\subset K$.
By this description,
we see that
$Q(K)\times_{Q(n)} Q(L)$ is isomorphic to $Q(K\cap L)$
for subcomplexes $K$ and $L$ of $\Delta^n$.

%We notice that 
%there is an isomorphism
%of scaled simplicial sets
%$Q(A)\times_{Q(C)} Q(B)\cong Q(A\cap B)$
%for subcomlexes $A$ and $B$ of $C$.
%\end{remark}

By
Lemmas~\ref{lemma:general-reedy-Segal-space} and 
\ref{lemma:twr-Z-Segal-space},
it suffices to show that
$Q(\bullet)$ satisfies the co-Segal condition
in order to prove Proposition~\ref{prop:lifting-SS}.

\begin{lemma}\label{lemma:inner-scaled-anodyne-Q}
The inclusion map
$Q(\Lambda^n_i)\to Q(n)$
is scaled anodyne for any $0<i<n$.
\end{lemma}

\proof
Let $(K^n_i)_{\S}$
be a subcomplex of $Q(n)$
defined in \cite[Construction 2.7]{Garcia-Stern},
that is,
$(K^n_i)_{\S}=Q(\Lambda^n_i)\cup \Delta^n_{\sharp}
\cup \Delta^{n,{\rm op}}_{\sharp}$.
By
\cite[Lemmas~2.10.1 and 2.12.1]{Garcia-Stern},
the inclusion map $(K^n_i)_{\S}\to Q(n)$ is scaled anodyne,
and hence it suffices to show that
$Q(\Lambda^n_i)\to (K^n_i)_{\S}$
is scaled anodyne.
This follows by observing that 
$Q(\Lambda^n_i)\cap (\Delta^n_{\sharp}\cup\Delta^{n,{\rm op}}_{\sharp})=
\Lambda^n_{i,\sharp}\cup \Lambda^{n,{\rm op}}_{i,\sharp}$
and
$\Lambda^n_{i,\sharp}\cup \Lambda^{n,{\rm op}}_{i,\sharp}
\to \Delta^n_{\sharp}\cup \Delta^{n,{\rm op}}_{\sharp}$
is scaled anodyne
by \cite[Remark~3.1.5]{Lurie3}.
\qed

\begin{lemma}\label{lemma:lambda-consecive-Q-anodyne}
For $0< i<n$,
the inclusion map
$Q(\Lambda^n_{\{1,\ldots, i\}})\to Q(n)$
is scaled anodyne.
\end{lemma}

\proof
We prove the lemma by induction on $n$.
When $n=2$,
the map $Q(\Lambda^2_1)\to Q(2)$ is scaled anodyne
by Lemma~\ref{lemma:inner-scaled-anodyne-Q}.
Now, we suppose that $n\ge 3$
and assume that the lemma holds for less than $n$.
We show that
$Q(\Lambda^n_{\{1,\ldots,i\}})\to Q(n)$
is scaled anodyne by induction of $i$.
When $i=1$,
it holds
by Lemma~\ref{lemma:inner-scaled-anodyne-Q}.
We suppose that $2\le i<n$ and
assume that it holds for less than $i$.
Since 
$\Lambda^{n}_{\{1,\ldots,i\}}\cap \Delta^{[n]-\{i\}}=
\Lambda^{[n]-\{i\}}_{\{1,\ldots,i-1\}}$,
there is a pushout diagram of scaled simplicial sets
\[ \xymatrix{
    Q(\Lambda^{[n]-{\{i\}}}_{\{1,\ldots,i-1\}})\ar[r]\ar[d]&
    Q(\Delta^{[n]-\{i\}})\ar[d]\\
    Q(\Lambda^n_{\{1,\ldots,i\}})\ar[r]&
    Q(\Lambda^n_{\{1,\ldots,i-1\}}).\\
    }\]
By using the isomorphism
$Q(\Delta^{[n]-\{i\}})\cong Q(n-1)$
and the hypothesis of induction,
the top horizontal arrow is
scaled anodyne, 
and hence so is the bottom horizontal arrow.
By the hypothesis of induction,
$Q(\Lambda^n_{\{1,\ldots,i-1\}})\to Q(n)$ is scaled anodyne.
Thus,
the composite
$Q(\Lambda^n_{\{1,\ldots,i\}})\to
Q(\Lambda^n_{\{1,\ldots,i-1\}})\to Q(n)$
is also scaled anodyne.
\qed

\proof[Proof of Proposition~\ref{prop:lifting-SS}]
It suffices to show that
$Q(\bullet)$ satisfies the co-Segal condition.
For $0\le i\le j\le n$,
we define
a subcomplex ${\rm sp}_{[i,j]}$ of $\Delta^n$ by
\[ {\rm sp}_{[i,j]}=\Delta^{\{i,i+1\}}
   \coprod_{\Delta^{\{i+1\}}}\cdots
   \coprod_{\Delta^{\{j-1\}}}
   \Delta^{\{j-1,j\}}. \]
%We write ${\rm sp}_n$ for
%${\rm sp}_{[0,n]}$ for simplicity.
We show that
the co-Segal map
$Q({\rm sp}_{[0,n]})\to Q(n)$
is scaled anodyne by induction on $n$.
When $n=0,1$,
it is trivial.
We suppose $n\ge 2$ and
assume that it holds for less than $n$.
We define a subcomplex $s(i,n)$ of $\Delta^n$ by 
\[ s(i,n)=\Delta^{\{0,\ldots,n-1\}}\coprod_{\Delta^{\{i,\ldots,n-1\}}}
   \Delta^{\{i,\ldots,n\}}\]
for $0\le i <n$.
By setting $S(i,n)=Q(s(i,n))$,
we obtain 
a filtration of scaled simplicial sets
$Q({\rm sp}_{[0,n]})\to S(n-1,n)\to S(n-2,n)\to \cdots \to
S(0,n)=Q(n)$.
The inclusion map
$Q({\rm sp}_{[0,n]})\to S(n-1,n)$
is scaled anodyne
since it is obtained
as a pushout 
of $Q({\rm sp}_{[0,n-1]})\to Q(\Delta^{\{0,\ldots,n-1\}})$,
which is scaled anodyne 
by the hypothesis of induction,
along the inclusion map
$Q({\rm sp}_{[0,n-1]})\to
Q({\rm sp}_{[0,n]})$.
%Q({\rm Sp}_{n-1})\coprod_{Q(\Delta^{\{n-1\}})}
%Q(\Delta^{\{n-1,n\}})$.
Therefore,
it suffices to show that
$S(i+1,n)\to S(i,n)$ is scaled anodyne
for $0\le i<n-1$. 
Since $s(i+1,n)\cap \Delta^{\{i,\ldots,n\}}=
\Lambda^{\{i,\ldots,n\}}_{\{i+1,\ldots,n-1\}}$,
there is a pushout diagram
\[ \xymatrix{
    Q(\Lambda^{\{i,\ldots,n\}}_{\{i+1,\ldots,n-1\}})
    \ar[r]\ar[d]&
    Q(\Delta^{\{i,\ldots,n\}})\ar[d]\\
    S(i+1,n)\ar[r] &  S(i,n).\\
}\]
The top horizontal arrow
is scaled anodyne
by Lemma~\ref{lemma:lambda-consecive-Q-anodyne},
and hence so is the bottom horizontal arrow.
\qed

\subsection{Completeness of the Segal space
$\btwr(Z)$}
\label{subsection:twrZ-complete-Q}

We recall that a Segal space $W$
is complete if the map
$s_0: W_0\to W_1^{\rm eq}$
is an equivalence,
where $W^{\rm eq}_1$ 
is the full subspace of $W_1$
spanned by equivalence morphisms.
The goal of this subsection is
to show that
the Segal space $\btwr(Z)$
is complete.

\begin{theorem}\label{thm:GHmodel-complete-Segal-TwZ}
For any $\infty$-bicategory $Z$,
the Segal space
$\btwr(Z)$ is complete.
\end{theorem}

First,
we recall relationships between
the model categories $\bsetdelta^{\rm CSS}$ and $\setjoy$,
where $\bsetdelta^{\rm CSS}$
is the category of bisimplicial sets $\bsetdelta$
equipped with the complete Segal space model structure.
Let $\pi_1: \Delta^{\rm op}\times\Delta^{\rm op}\to \Delta^{\rm op}$
be the first projection,
and let $i_1: \Delta^{\rm op}\cong \{[0]\}\times\Delta^{\rm op}
\to \Delta^{\rm op}\times\Delta^{\rm op}$ be the inclusion. 
These functors induce an adjunction
$\pi_1^*: \setdelta\rightleftarrows \bsetdelta: i_1^*$,
where $(\pi_1^*A)_{m,n}=A_m$ and
$(i_1^*B)_{n}=B_{n,0}$
for $A\in \setdelta$ and $B\in\bsetdelta$.
By
\cite[Theorem~4.11]{Joyal-Tierney},
the adjunction
$(\pi_1^*,i_1^*)$ forms 
a Quillen equivalence
\[ \pi_1^*: \setjoy\rightleftarrows \bsetdelta^{\rm CSS}: i_1^*.\]

There is a natural isomorphism of sets
$\mapsc(X,Y)_0\cong
{\rm Hom}_{\setsc}(X,Y)$
for a scaled simplicial set $X$ and
an $\infty$-bicategory $Y$.
This implies that there is a natural isomorphism
of simplicial sets
\[ i_1^*\btwr(Z)\cong \twr(Z) \]
for any $\infty$-bicategory $Z$.

We denote by $E(n)$
the nerve of
the groupoid
$[n]^{\rm gpd}$
freely generated by the category $[n]$.
Let $t: \Delta\times\Delta\to \setdelta$
be the functor given by
$t([m],[n])=\Delta^m\times E(n)$.
By a left Kan extension of $t$
along the Yoneda embedding $\Delta\times\Delta\to \bsetdelta$,
we obtain a functor $t_!: \bsetdelta\to \setdelta$.
The functor $t_!$ admits a right adjoint
$t^!: \setdelta\to \bsetdelta$ given by
$t^!(A)_{m,n}={\rm Hom}_{\setdelta}
(\Delta^m\times E(n), A)$
for $A\in\setdelta$.
By \cite[Theorem~4.12]{Joyal-Tierney},
the adjunction $(t_!,t^!)$ forms
another Quillen equivalence 
\[ t_!: \bsetdelta^{\rm CSS}\rightleftarrows \setjoy: t^!. \]

Next,
we will construct a lifting of the functor
$u_1: \setsc\to \setjoy$
to $\bsetdelta^{\rm CSS}$.
We define a functor
$v_1: \setsc\to\bsetdelta^{\rm CSS}$
to be the composite $t^!u_1$.
Since $t^!$ and $u_1$ are right Quillen functors,
so is $v_1$.
Since there is a natural isomorphism
$i_1^*t^!A\cong A$ for any simplicial set $A$,
$v_1$ is a lifting of $u_1$
through $i_1^*$ in the sense that
there is a natural isomorphism of functors
$i_1^*v_1\cong u_1$.

For a scaled simplicial set $X$,
we define a bisimplicial set
$v_1^{\rm rev}X$ by
$v_1^{\rm rev}X= t^! (u_1X^{\rm op})$.
We notice that 
there is a natural isomorphism of sets
$(v_1X)_{m,n}\cong 
\mapsc(\Delta^m_{\sharp},X)_n$ 
and 
$(v_1^{\rm rev}X)_{m,n}\cong
\mapsc(\Delta^{m,{\rm op}}_{\sharp},X)_n$
for any $\infty$-bicategory $X$.

For an $\infty$-bicategory $Z$,
the inclusion maps
$\Delta^n_{\sharp}\coprod\Delta^{n,{\rm op}}_{\sharp}
\to   Q(n)$
for $n\ge 0$
induce a map of Segal spaces
\[ q: \btwr(Z)\longrightarrow
   v_1Z\times v_1^{\rm rev}Z.\]
The map $q$ induces a commutative diagram
\begin{equation}\label{eq:complete-edge-diagram-Q}
   \xymatrix{
    \btwr(Z)_0\ar[rr]^{s_0}\ar[d]_{q_0}&&
    \btwr(Z)_1^{\rm eq}\ar[d]^{q_1^{\rm eq}}\\
    (v_1Z)_0 \times (v_1^{\rm rev}Z)_0
    \ar[rr]^{s_0\times s_0}&&
    (v_1Z)_1^{\rm eq} \times (v_1^{\rm rev}Z_1)^{\rm eq}.\\
}
\end{equation}
We notice that
the bottom horizontal arrow
is an equivalence
since the Segal spaces $v_1Z$ and
$v_1^{\rm rev}Z$ are complete. 
Hence,
in order to prove that
the top horizontal arrow is an equivalence,
it suffices to show that
(\ref{eq:complete-edge-diagram-Q})
is a pullback diagram
in the $\infty$-category of $\infty$-groupoids.

%In order to prove Theoreom~\ref{thm:GHmodel-complete-Segal-TwZ},
%we make preliminary results.
We take a $0$-simplex $M$ in the Kan complex
$\btwr(Z)_1=\mapsc(Q(1),Z)$
that is represented
by a map $M: Q(1)\to Z$
of scaled simplicial sets.
We put $N_0=d_1(M)$ and $N_1=d_0(M)$.
We would like to have a necessary and sufficient
condition for
$M$ being an equivalence
in the Segal space $\btwr(Z)$.

\begin{proposition}\label{prop:twr-equivalence-morphism-Q}
A morphism $M\in \btwr(Z)_1$
is an equivalence
in the Segal space $\btwr(Z)$ 
if and only if
$M$ extends to a map from 
$(\Delta^1\star\Delta^{1,{\rm op}})_{\sharp}$ and
the edges $M(\Delta^{\{00,01\}})$ and
$M(\Delta^{\{11,10\}})$
are equivalence $1$-morphisms in $Z$.
\end{proposition}

First,
we show the ``if'' direction
of Proposition~\ref{prop:twr-equivalence-morphism-Q}.

\begin{lemma}\label{lemma:M-eq-if-direction-Q}
If $M\in\btwr(Z)_1$ is an equivalence,
then $M$ extends to a map from 
$(\Delta^1\star\Delta^{1,{\rm op}})_{\sharp}$ and
the edges $M(\Delta^{\{00,01\}})$ and
$M(\Delta^{\{11,10\}})$
are equivalence $1$-morphisms in $Z$.
\end{lemma}

\proof
Suppose that $M\in \btwr(Z)_1^{\rm eq}$.
%with $d_0(M)=N_1$ and $d_1(M)=N_0$.
%By the right vertical arrow
%in (\ref{eq:complete-edge-diagram}),
By applying
  $q_1^{\rm eq}$
to $M$,
we see that 
$M(\Delta^{\{00,01\}})$ and $M(\Delta^{\{11,10\}})$
are equivalence $1$-morphisms in $Z$.

We shall show that
$M$ extends to a map from
$(\Delta^1\star\Delta^{1,{\rm op}})_{\sharp}$.
For this purpose,
it suffices to show that
$M(\Delta^{\{00,10,11\}})$ is thin
by \cite[Remark~3.1.4]{Lurie3}.
By the assumption,
there are $L,R\in \btwr(Z)_2$ such that
$d_0(L)=M, d_1(L)=s_0(N_1)$ and
$d_1(R)=s_0(N_0),d_2(R)=M$.
%By applying Lemma~\ref{lemma:inner-twsc},
By using Lemma~\ref{lemma:inner-scaled-anodyne-Q}
for $n=3$ and $i=2$,
we obtain $A\in \btwr(Z)_3$
such that $d_0(A)=R$,
$d_3(A)=L$,
which is given by a map
$A: Q(3)\to Z$
of scaled simplicial sets.
%Since $A(\Delta^{\{00,01\}})$
%is an equivalence,
%$M(\Delta^{\{00,10,11\}})$ is thin
%if and only if
%$A(\Delta^{\{00,11,12\}})$ is thin
%by \cite[Remark~3.1.4]{Lurie3} and
%\cite[Proposition~3.4(1)]{GHL}.
We notice that
$A(\Delta^{\{00,10,12\}})$ and
$A(\Delta^{\{01,11,13\}})$ are thin triangles
since $d_1d_3(A)=d_1(L)=s_0(N_1)$
and $d_1d_0(A)=d_1(R)=s_0(N_0)$
are degenerate.
Furthermore,
we see that
the edges $A(\Delta^{\{00,01\}})$
and $A(\Delta^{\{10,11\}})$
are equivalences
by applying $q_3$ to $A$. 

We would like to show that 
$M(\Delta^{\{00,10,11\}})=A(\Delta^{\{01,11,12\}})$
is thin.
We consider the restriction of $A$
to the $3$-simplex
$\Delta^{\{00,01,11,12\}}$.
In order to show that
$A(\Delta^{\{01,11,12\}})$ is thin,
it suffices to show that 
$A(\Delta^{\{00,11,12\}})$ is thin
by \cite[Proposition~3.4(1)]{GHL}
since
the edge $A(\Delta^{\{00,01\}})$
is an equivalence and
the triangle
$A(\Delta^{\{00,01,11\}})$
is thin.

First,
we will show that $A(\Delta^{\{00,11,13\}})$ is thin.
We consider the restriction of $A$
to the $3$-simplex $\Delta^{\{00,01,11,13\}}$.
Since $A(\Delta^{\{01,11,13\}}),
A(\Delta^{\{00,01,13\}}), A(\Delta^{\{00,01,11\}})$
are thin,
so is $A(\Delta^{\{00,11,13\}})$
by \cite[Remark~3.1.4]{Lurie3}.
Next,
we consider the restriction
of $A$ to the $4$-simplex
$\Delta^{\{00,10,11,12,13\}}$.
We notice that
the trivial cofibrations
are invariant
under taking opposites of simplicial sets
since the bicategorical equivalences are
invariant (cf.~\cite[Remark~1.32]{GHL}).
A lifting against
the opposite of An2,
which is a trivial cofibration,
implies that $A(\Delta^{\{00,10,11\}})$ is thin.
Finally,
we consider the restriction of $A$
to the $3$-simplex $\Delta^{\{00,10,11,12\}}$.
Since
  $A(\Delta^{\{00,10,12\}}),
  A(\Delta^{\{00,10,11\}}),
  A(\Delta^{\{10,11,12\}})$
  are thin and
the edge $A(\Delta^{\{10,11\}})$ is
an equivalence,
$A(\Delta^{\{00,11,12\}})$ is thin
by \cite[Proposition~3.4(2)]{GHL}.
\qed

\bigskip

Next,
we show
the ``only if'' direction
of Proposition~\ref{prop:twr-equivalence-morphism-Q}.

\begin{lemma}\label{lemma:M-eq-only-if-direction-Q}
If $M\in\btwr(Z)_1$ extends to a map from 
$(\Delta^1\star\Delta^{1,{\rm op}})_{\sharp}$ and
the edges $M(\Delta^{\{00,01\}})$ and
$M(\Delta^{\{11,10\}})$
are equivalence $1$-morphisms in $Z$,
then $M$ is an equivalence.
\end{lemma}

\proof
It suffices to construct $L,R\in \btwr(Z)_2$
such that
$d_0(L)=M, d_1(L)=s_0(N_1)$ and
$d_1(R)=s_0(N_0), d_2(R)=M$.
We recall that
we have a cartesian fibration
$p: \twr(Z)\to u_1Z\times u_1Z^{\rm op}$
of $\infty$-categories
by \cite[Theorem~0.1]{Garcia-Stern}.
Using the fact that the set of $n$-simplices
of $\twr(Z)$
is ${\rm Hom}_{\setsc}(Q(n),Z)$,
we regard a $0$-simplex $X\in \btwr(Z)_n=\mapsc(Q(n),Z)$
as an $n$-simplex $\overline{X}\in \twr(Z)_n$.
In particular,
we regard $M\in \btwr(Z)_1=\mapsc(Q(1),Z)$
as a morphism $\overline{M}$ of $\twr(Z)$.
By \cite[Theorem~2.6]{Garcia-Stern},
a morphism in $\twr(Z)$
is $p$-cartesian if and only if
it extends to a map from
$(\Delta^1\star\Delta^{1,{\rm op}})_{\sharp}$.
Thus,
$\overline{M}$ is 
a $p$-cartesian morphism of $\twr(Z)$.
The assumption that
$M(\Delta^{\{00,01\}})$ and
$M(\Delta^{\{11,10\}})$
are equivalence $1$-morphisms in $Z$
implies that
$\overline{M}$ is an equivalence in $\twr(Z)$
by \cite[Proposition~2.4.1.5]{Lurie1}.
Hence there exist 
$\overline{L},\overline{R}\in \twr(Z)_2$
such that
$d_0(\overline{L})=\overline{M},
d_1(\overline{L})=s_0(\overline{N}_1)$ and
$d_1(\overline{R})=s_0(\overline{N}_0),
d_2(\overline{R})=\overline{M}$.
By using the isomorphism of simplicial sets
$i_1^*\btwr(Z)\cong \twr(Z)$,
we obtain the desired 
$L,R\in \btwr(Z)_2$.
\qed
%\bigskip

\proof[Proof of Proposition~\ref{prop:twr-equivalence-morphism-Q}]
%Proposition~\ref{prop:twr-equivalence-morphism}
The proposition follows from
Lemmas~\ref{lemma:M-eq-if-direction-Q}
and \ref{lemma:M-eq-only-if-direction-Q}.
%Lemmas~\ref{lemma:M-eq-if-direction},
%\ref{lemma:M-eq-only-if-const-L}
%and \ref{lemma:M-eq-only-if-const-R}.
\qed

\bigskip

Now,
we turn to the proof of
Theorem~\ref{thm:GHmodel-complete-Segal-TwZ}.
%Theorem~\ref{thm:CSS-complete-Segal-TwZ}.
For this purpose,
we give some preliminary lemmas.
%We recall that
%$\sigma(k)=\Delta^{\{00,\ldots,0k,1k,\ldots,nk\}}$
%is a simplex of
%$\Omega^n$ for $0\le k\le n$.
%
%We set $X=\Delta^{\{00,01\}}\cup\Delta^{\{01,11\}}\cup
%\Delta^{\{10,11\}}$.
%Since the inclusion map
%$X\hookrightarrow \Omega^1_{\sharp}$
%is obtained by iterated pushouts
%along the scaled anodyne maps {\bf An1},
%we obtain the following lemma.
%
We set
\[ \begin{array}{rcl}
    P(1)&=&(\Delta^1\star\Delta^{1,{\rm op}})_{\sharp},\\[2mm]
     {\rm Sp}_1(P)&=&\Delta^{\{00,01\}}_{\sharp}\coprod_{\Delta^{\{01\}}_{\sharp}}
     \Delta^{\{01,11\}}_{\sharp}\coprod_{\Delta^{\{11\}}_{\sharp}}
     \Delta^{\{11,10\}}_{\sharp}.\\[2mm]
   \end{array}\]

\begin{lemma}\label{lemma:scaled-anodyne-X-omega-1-sharp-Q}
The inclusion map
${\rm Sp}_1(P)\to P(1)$
%$\Delta^{\{00,01\}}_{\sharp}\coprod_{\Delta^{\{01\}}_{\sharp}}
%\Delta^{\{01,11\}}_{\sharp}\coprod_{\Delta^{\{11\}}_{\sharp}}
%\Delta^{\{11,10\}}_{\sharp}
%\to P(1)$
is scaled anodyne.
\end{lemma}

\proof
This follows from the fact that
it is obtained by iterated pushouts
along scaled anodyne maps of type An1.
\qed

\bigskip

We set
\[ \widetilde{P}(1)=
   \Delta^0_{\sharp}\coprod_{\Delta^{\{00,01\}}_{\sharp}}
   P(1)
   \coprod_{\Delta^{\{10,11\}}_{\sharp}}\Delta^0_{\sharp}.\]
The map $s^0: Q(1)\to Q(0)$ extends
to a map $\widetilde{s}^0:\widetilde{P}(1)\to Q(0)$
of scaled simplicial sets.
We show that
it is a bicategorical equivalence.

\begin{lemma}\label{lemma:s-zero-sharp-bicat-eq-Q}
The map $\widetilde{s}^0:\widetilde{P}(1)\to Q(0)$
is a bicategorical equivalence.
\end{lemma}

\proof
There are canonical maps
$Q(1)\to P(1)\to \widetilde{P}(1)$.
We denote by $\widetilde{d}{}^0: Q(0)\to \widetilde{P}(1)$
the composite of these maps
with $d^0: Q(0)\to Q(1)$.
In order to prove the lemma,
it suffices to show that
$\widetilde{d}^0$ is a trivial cofibration
since the composite $\widetilde{s}^0\widetilde{d}^0$
is the identity of $Q(0)$.

We set
\[ {\rm Sp}_1(\widetilde{P})=
     \Delta^0_{\sharp}\coprod_{\Delta^{\{00,01\}}_{\sharp}}
     {\rm Sp}_1(P)
     \coprod_{\Delta^{\{10,11\}}_{\sharp}}\Delta^0_{\sharp}.\]
%We will show that $\widetilde{d}^0$ is a trivial cofibration.
%By Lemma~\ref{lemma:d--sharp-trivial-cof},
By Lemma~\ref{lemma:scaled-anodyne-X-omega-1-sharp-Q},
we obtain a scaled anodyne map
${\rm Sp}_1(\widetilde{P})\to \widetilde{P}(1)$.
The desired result follows from
the fact that the map $d^0$
induces an isomorphism
$Q(0)\cong {\rm Sp}_1(\widetilde{P})$
of scaled simplicial sets.
\qed

\proof[Proof of Theorem~\ref{thm:GHmodel-complete-Segal-TwZ}]
We would like to show that
(\ref{eq:complete-edge-diagram-Q})
is a pullback diagram
in the $\infty$-category of
$\infty$-groupoids.
First,
we note that
the map
$q_n: \btwr(Z)_n\to (v_1Z)_n\times (v_1^{\rm rev}Z)_n$
is a Kan fibration between Kan complexes
for any $n\ge 0$
by Proposition~\ref{prop:cof-induce-map-sc-fibration}.
By using the fact that
$W_1^{\rm eq}$ is a full subspace
of the Kan complex $W_1$
for a Segal space $W$,
we see that
$q_1^{\rm eq}: \btwr(Z)_1^{\rm eq}\to
(v_1Z)_1^{\rm eq}\times (v_1^{\rm rev}Z)_1^{\rm eq}$
is also a Kan fibration between Kan complexes.
Thus,
in order to prove that
(\ref{eq:complete-edge-diagram-Q})
is a pullback diagram,
it suffices to show that 
the induced map
$\btwr(Z)_{0,(x,y)}\to \btwr(Z)_{1,(s_0(x),s_0(y))}^{\rm eq}$
on fibers
is an equivalence
for any $(x,y)\in (v_1Z)_0\times (v_1^{\rm rev}Z)_0$.

%We set $W=     \Delta^0\coprod_{\Delta^{\{00,01\}}}
%     \Omega^1_{\sharp}
%           \coprod_{\Delta^{\{10,11\}}}\Delta^0$.
The map
$\widetilde{s}^0: \widetilde{P}(1)\to Q(0)$
induces a map 
$s_0: \btwr(Z)_0\to \mapsc(\widetilde{P}(1),Z)$,
which is an equivalence
by Lemma~\ref{lemma:s-zero-sharp-bicat-eq-Q}.
By composition with the inclusion map
$\Delta^0_{\sharp}\coprod\Delta^0_{\sharp}
\to \widetilde{P}(1)$,
we obtain
a Kan fibration between Kan complexes
$\widetilde{q}_1:
\mapsc(\widetilde{P}(1),Z)
\to (v_1Z)_0\times (v_1^{\rm rev}Z)_0$
by Proposition~\ref{prop:cof-induce-map-sc-fibration},
which makes the following diagram commute
\[ \xymatrix{
   \btwr(Z)_0\ar[rr]^{s_0}\ar[dr]_{p_0} &&
   \mapsc(\widetilde{P}(1),Z)
        \ar[dl]^{\widetilde{q}_1}\\
   &(v_1Z)_0\times (v_1^{\rm rev}Z)_0.&}\]
By taking fibers over $(x,y)\in
(v_1Z)_0\times (v_1^{\rm rev}Z)_0$,
we obtain the desired equivalence.
\qed

%%%\newpage
%\input{app-simple-twist}
\section{Another model of
twisted arrow $\infty$-categories}
\label{section:model-twisted-arrow}

In this section
we give another simple model
of twisted arrow $\infty$-categories
for an $\infty$-bicategory.
In \S\ref{subsection:introduction-twscbullet}
we introduce a cosimplicial scaled simplicial set
$\twscbullet$,
where the scaled simplicial set $T(n)$ is a subcomplex
of $Q(n)$ for all $n\ge 0$.
In \S\ref{subsection:SS-twrZ}
we define a bisimplicial set
$\bbtwr(Z)$ to be
$\{\mapsc(T(n),Z)\}_{n\ge 0}$
for an $\infty$-bicategory $Z$.
We show that the bisimplicial set
$\bbtwr(Z)$ is a Segal space. 
In \S\ref{subsection:Segal-sp-twrZ-complete}
we show that
the cosimplicial objects $T(\bullet)$
and $Q(\bullet)$ are equivalent
in the Reedy model structure
on the category of cosimplicial scaled simplicial sets.
By using this result,
we show 
that the Segal space $\bbtwr(Z)$ is complete.
By applying the right Quillen equivalence
$i_1^*: \bsetdelta^{\rm CSS}\to \setjoy$
to $\bbtwr(Z)$,
we obtain an $\infty$-category
$\bstwr(Z)$ which is equipped with a map
$p_{\mathbb{T}}: \bstwr(Z)\to u_1Z\times u_1Z^{\rm op}$.
In \S\ref{subsection:cartesian-fibration-structure}
we show that $p_{\mathbb{T}}$ is a cartesian fibration
equivalent to $p: \twr(Z)\to u_1Z\times u_1Z^{\rm op}$.
In particular,
$p_{\mathbb{T}}$ is classified by the restricted
mapping $\infty$-category functor.
%In \S\ref{subsection:Segal-sp-twrZ-complete}
%by showing that $\rtwr(Z)_{\bullet}$ and $\twr(Z)_{\bullet}$
%are equivalent,
%we show that 
%$\twr(Z)_{\bullet}$ is a complete Segal complete.
%In \S\ref{subsection:model-twisted-arrow}
%we show that
%$\rtwr(Z)$ equipped with the cartesian fibration $\mathbf{p}$
%is equivalent to $\twr(Z)$ with $p$
%and hence $\mathbf{p}$
%is classified by the restricted
%mapping $\infty$-category functor.

\subsection{The cosimplicial scaled simplicial set
$\twscbullet$}
\label{subsection:introduction-twscbullet}

In this subsection
we introduce a cosimplicial object
$\twscbullet$ of scaled simplicial sets.
%For an $\infty$-bicategory $Z=(\overline{Z},T_Z)$,
%we construct a complete Segal space
%$\rtwr(Z)_{\bullet}$ that
%is equivalent to
%$\twr(Z)_{\bullet}$.

\begin{definition}\label{definitoin:omega-n}
\rm
For $n\ge 0$,
we consider the simplicial set
$\Delta^n\star\Delta^{n,{\rm op}}$.
We recall that we represent it as the following diagram
%\if0
\[ \xymatrix{
  00\ar[r]\ar[d]
  & 01\ar[r]\ar[d]
  & \cdots\ar[r]
  & 0n\ar[d] \\
  10
  & 11\ar[l]
  & \cdots\ar[l]
  & 1n.\ar[l]\\
}   \]
%\fi
%For a subcomplex $K$ of $\Delta^n\star\Delta^{n,{\rm op}}$,
%we denote by $K_{\dagger}$
%the scaled simplicial set whose underlying simplicial set
%is $K$ equipped with the induced scaling from $Q(n)$.

For a vertex
$v=(i,r)$ of $\Delta^1\times \Delta^n$,
we set $\hat{v}=ir$.
For a simplex $\sigma$ of $\Delta^1\times\Delta^n$,
we denote by $\hat{\sigma}$
the simplex of $\Delta^n\star\Delta^{n,{\rm op}}$
spanned by $\hat{v}$
for vertices $v$ of $\sigma$.
For a subcomplex $K$ of $\Delta^1\times \Delta^n$,
we denote by $\Omega(K)$
the subcomplex of $\Delta^n\star\Delta^{n,{\rm op}}$
spanned by $\hat{\sigma}$
for $\sigma\in K$.
We write $\Omega^n$
for $\Omega(\Delta^1\times\Delta^n)$
for simplicity.
We notice that
$\Omega^n$ is an $(n+1)$-dimensional subcomplex
of $\Delta^n\star\Delta^{n,{\rm op}}$.
The collection
$\Omega^{\bullet}=\{\Omega^n\}_{n\ge 0}$
forms a cosimplicial object of simplicial sets
by restricting the cosimplicial structure on
$\{\Delta^n\star\Delta^{n,{\rm op}}\}_{n\ge 0}$.

We define a scaled simplicial set $T(n)$.
The underlying simplicial set
of $T(n)$ is $\Omega^n$.
In other words,
it is a full subcomplex of $\Delta^n\star\Delta^{n,{\rm op}}$
spanned by
$(n+1)$-dimensional simplices
\[ \sigma(r)=\Delta^{\{00,\ldots,0r,1r,\ldots,1n\}} \]
for $0\le r\le n$.
The scaling of $T(n)$ is induced by $Q(n)$.
Concretely,
the set $T_{T(n)}$
of thin triangles is given by
\[ \begin{array}{rccl}
  T_{T(n)}
    &=& &
    \{\Delta^{\{ik,ik',ik''\}}|\
    i=0,1,\
    0\le k< k'< k'' \le n\}\\[2mm]
    & & \cup &
    \{\Delta^{\{0k,0k',1k''\}}|\
    0\le k<  k'\le k'' \le n\}\\[2mm]
  & & \cup &
    \{\mbox{\rm degenerate}\}.     
\end{array}\]
We note that
$(\Omega^n)_2\cap T3=\emptyset$
since $\sigma(r)$ does not contain any
triangles in $T3$
for $0\le r\le n$,
where $T3=\{\Delta^{\{1k,1k',0k''\}}|\ 0\le k<k'\le k''\le n\}$
is a subset of $T_{Q(n)}$.

We can verify that
the collection $\twscbullet=\{T(n)\}_{n \ge 0}$
forms a cosimplicial object
in the category of scaled simplicial sets
by restricting the cosimplicial structure on $Q(\bullet)$.
\end{definition}

\subsection{The Segal space
$\bbtwr(Z)$}
\label{subsection:SS-twrZ}

In this subsection
we introduce a bisimplicial set
$\bbtwr(Z)$
for an $\infty$-bicategory $Z$,
and show that it is a Segal space.
We will show that
it is a complete Segal space
in \S\ref{subsection:Segal-sp-twrZ-complete} below.
%We will show that
%it is a model of twisted arrow $\infty$-category
%in \S\ref{subsection:cartesian-fibration-structure} below.
%%%in \S\ref{subsection:model-twisted-arrow}.

\begin{definition}\rm
Let $Z$ be an $\infty$-bicategory.
We set
\[ \bbtwr(Z)_n=\mapsc(T(n),Z). \]
%where $\mapsc(-,-)$ is the mapping space
%on $(\setsc)^{\circ}$.
The collection
$\bbtwr(Z)=\{\bbtwr(Z)_n\}_{n\ge 0}$
forms a bisimplicial set
by using the cosimplicial structure on
$T(\bullet)$
\end{definition}

The goal of this subsection is
to show the following proposition.

\begin{proposition}\label{prop:twr-complete-Segal}
For any $\infty$-bicategory $Z$,
the bisimplicial set
$\bbtwr(Z)$ is a Segal space.
\end{proposition}

In order to prove Proposition~\ref{prop:twr-complete-Segal},
we give some preliminary lemmas.
%We consider the Reedy model structure
%on the category of cosimplicial objects
%in $\setsc$
First,
we show that $\twscbullet$
is Reedy cofibrant.

\begin{lemma}\label{lemma:twsc-reedy-cofibrant}
The cosimplicial object
$\twscbullet$ in $\setsc$
is Reedy cofibrant.
\end{lemma}

\proof
This follows by observing
that the $n$th latching object
is isomorphic to
$\Omega(\Delta^1\times\partial\Delta^n)_{\S}$.
\qed

\bigskip

For $0\le i\le n$,
we define subcomplexes
$\Lambda^n_iT$ and
$\overline{\Lambda}{}^n_iT$
of $T(n)$ by
\[ \begin{array}{rcl}
   \Lambda^n_iT &=&
   \Omega(\Delta^1\times\Lambda^n_i)_{\S},\\[2mm]
   \overline{\Lambda}{}^n_iT &=&
   \Omega(\Delta^1\times\Lambda^n_i\cup
   \partial\Delta^1\times\Delta^n)_{\S}.\\
   \end{array}\]   

\begin{lemma}\label{lemma:inner-twsc}
For any $0<i<n$,
the inclusion map
$\Lambda^n_iT\to T(n)$
is scaled anodyne.
\end{lemma}

\proof
First,
we show that
the inclusion map
$\Lambda^n_iT\to
\overline{\Lambda}{}^n_iT$
is scaled anodyne.
There is a pushout diagram
\[ \xymatrix{
    \Omega(\partial\Delta^1\times \Lambda^n_i)_{\S}
    \ar[r]\ar[d]&
    \Omega(\partial\Delta^1\times\Delta^n)_{\S}\ar[d]\\
    \Lambda^n_iT
    \ar[r] & \overline{\Lambda}^n_iT}\]
of scaled simplicial sets.
We note that
there is a natural isomorphism
of scaled simplicial sets
$\Omega(\partial\Delta^1\times K)_{\S}\cong
K_{\sharp}\coprod K_{\sharp}^{\rm op}$
for any subcomplex $K$ of $\Delta^n$.
Since the top horizontal arrow is scaled anodyne
by \cite[Remark~3.1.5]{Lurie3},
so is the bottom horizontal arrow.
%This follows from the fact that
%it is obtained by taking iterated pushouts along
%the map $\Lambda^n_{i,\sharp}\hookrightarrow
%\Delta^n_{\sharp}$.

Next,
we show that
the inclusion map
$\overline{\Lambda}{}^n_iT
\to T(n)$
is scaled anodyne.
For $0\le r\le n$,
we recall that $\sigma(r)$ is
the $(n+1)$-dimensional simplex
$\Delta^{\{00,\ldots,0r,1r,\ldots,1n\}}$
of $\Delta^n\star\Delta^{n,{\rm op}}$.
For $0\le s\le n$,
we define a scaled simplicial set $X(s)$,
which is a subcomplex of $\Omega^n$
equipped with the induced scaling from $T(n)$, 
by
\[ X(s)= \overline{\Lambda}{}^n_iT
         \cup  \bigcup_{s\le r\le n}\
         \sigma(r)_{\S}.\]
By definition,
we obtain a filtration
\[ \overline{\Lambda}{}^n_iT=
   X(n+1)\to X(n)\to
   X(n-1)\to\cdots
   \to X(0) = T(n)\]
of scaled simplicial sets,
which consists of subcomplexes of $T(n)$.
In order to prove that
$\overline{\Lambda}{}^n_iT
\to T(n)$
is scaled anodyne,
it suffices to show that the inclusion map
$X(s+1)\to X(s)$
is scaled anodyne for any $0\le s\le n$.

We fix $s$ with $0\le s\le n$
and show that $X(s+1)\to X(s)$ is scaled anodyne.
We set $I(s)=\{00,\ldots,0s,1s,\ldots,1n\}$.
Then  we have $\sigma(s)=\Delta^{I(s)}$.
The intersection of
$\sigma(s)$ with
the underlying simplicial set of $X(s+1)$
has the form
$\Lambda^{I(s)}_{M(s)}$,
where
\[ M(s)=\left\{\begin{array}{ll}
                \{1i\} & (s=0),\\
                \{0s,1i\} & (0<s<i),\\
                \{0i\} & (s=i),\\
                \{0i,0s\} & (i<s\le n).\\
               \end{array} \right. \]
%We identify $\sigma(s)$
%with a map $\sigma(s):
%\Delta^{n+1}_{\S}\to X(s)$,
%where $\Delta^{n+1}_{\S}$ is $\Delta^{n+1}$
%equipped with the induced scaling.
Hence,
there is a pushout diagram
\[ \xymatrix{
  \Lambda^{I(s)}_{M(s),\S}\ar[r]\ar[d]&
    \Delta^{I(s)}_{\S}\ar[d]\\%[1mm]
     X(s+1)\ar[r] & X(s).\\
   }\]  
Thus,
it suffices to show that
the inclusion map
$\Lambda^{I(s)}_{M(s),\S}
\to \Delta^{I(s)}_{\S}$
is scaled anodyne. 
%For $0<s\le n$,
This follows from
Lemma~\ref{lemma:simplex-scaled-anodyne}.
%Lemmas~\ref{lemma:simplex-scaled-anodyne} and
%\ref{lemma:simplex-scaled-anodyne-2}.
\qed

\bigskip

Next,
we show that
$\twscbullet$ satisfies the co-Segal condition.

\begin{lemma}\label{lemma:twsc-co-Segal}
The cosimplicial scaled simplicial set
$\twscbullet$ satisfies the co-Segal condition.
\end{lemma}

\proof
It suffices to show that
the co-Segal map
\[ T(1)\coprod_{T(0)}\cdots
   \coprod_{T(0)}T(1)\longrightarrow T(n) \]
is scaled anodyne for all $n\ge 0$.
This follows by
induction on $n$
together with Lemma~\ref{lemma:inner-twsc}.
\qed

\begin{proposition}\label{prop:rtwrZ-Segal-space}
For any $\infty$-bicategory $Z$,
the bisimplicial set
$\bbtwr(Z)$ is a Segal space.
\end{proposition}

\proof
This follows from
Lemmas~\ref{lemma:general-reedy-Segal-space},
\ref{lemma:twsc-reedy-cofibrant},
and \ref{lemma:twsc-co-Segal}.
\qed

\subsection{Completeness of
the Segal space $\bbtwr(Z)$}
\label{subsection:Segal-sp-twrZ-complete}

The goal of this subsection
is to show that
the Segal space
$\bbtwr(Z)$ is complete.
For this purpose,
we prove that 
$T(\bullet)$ and $Q(\bullet)$
are equivalent in the Reedy model structure
on the category of cosimplicial scaled simplicial sets.

\begin{proposition}\label{prop:proof-of-equivalnece-model}
The inclusion map $\twscbullet\to Q(\bullet)$
is a weak equivalence in the Reedy model structure
on the category of cosimplicial objects in $\setsc$,
that is,
it is a levelwise bicategorical equivalence.
\end{proposition}

Using Proposition~\ref{prop:proof-of-equivalnece-model},
we can show that
the Segal space $\bbtwr(Z)$
is complete.

\begin{theorem}
\label{thm:twrZ-complete-Segal-space}
For any $\infty$-bicategory $Z$,
$\bbtwr(Z)$ is a complete Segal space.
\end{theorem}

\proof
%[Proof of Theorem~\ref{thm:twrZ-complete-Segal-space}]
The map
$\twscbullet\to Q(\bullet)$
of cosimplicial scaled simplicial sets
induces a map of Segal spaces
$\btwr(Z)  \to\bbtwr(Z)$,
which is a levelwise weak equivalence
by Proposition~\ref{prop:proof-of-equivalnece-model}.
The theorem follows from
the fact that $\btwr(Z)$
is a complete Segal space
by Theorem~\ref{thm:GHmodel-complete-Segal-TwZ}.
\qed

\bigskip

Now,
we turn to the proof of
Proposition~\ref{prop:proof-of-equivalnece-model}.
%For this purpose,
%we make some preliminary results.
For $0\le k<l\le n+1$,
we denote by $\tau(k,l)$
the simplex of $\Delta^n\star\Delta^{n,{\rm op}}$ given by
\[ \tau(k,l)=\Delta^{\{00,\ldots, 0k,0l,\ldots,0n,
                      1k,\ldots,1n\}}. \]
We understand $\tau(k,n+1)=
\Delta^{\{00,\ldots,0k,1k,\ldots,1n\}}$.
%We denote by $\tau(k,l)_{\dagger}$
%the scaled simplicial set
%whose underlying simplicial set is $\tau(k,l)$
%equipped with the induced scaling from $Q(n)$. 

\proof[Proof of
Proposition~\ref{prop:proof-of-equivalnece-model}]
For $0\le k\le n$,
we define a scaled simplicial set
$U(k)$,
which is a subcomplex of $\Delta^n\star\Delta^{n,{\rm op}}$
equipped with the induced scaling from $Q(n)$,
by
\[ U(k)=T(n)\cup \bigcup_{k\le i\le n}
   \tau(i,i+1)_{\S}. \]
By definition,
we obtain a filtration
$T(n)= U(n)\to\cdots\to
%X(1)=\twscn\cup \tau(1,2)_{\dagger}$.
U(0)=Q(n)$
of scaled simplicial sets
consisting of subcomplexes of $Q(n)$.
In order to prove
Proposition~\ref{prop:proof-of-equivalnece-model},
it suffices to show that
$U(k+1)\to U(k)$ is scaled anodyne
for all $0\le k<n$.

For this purpose,
we introduce a further filtration.
For $k+1\le l\le n+1$,
we define a scaled simplicial set
$V(k,l)$,
which is a subcomplex of $\Delta^n\star\Delta^{n,{\rm op}}$
equipped with the induced scaling
from $Q(n)$, by
\[ V(k,l)=U(k+1)\cup \bigcup_{l\le j\le n+1}\tau(k,j)_{\S}. \]
By definition,
we obtain a filtration
$U(k+1)=V(k,n+1)\to\cdots\to V(k,k+1)=U(k)$
of scaled simplicial sets
consisting of subcomplexes of $U(k)$.
In order to prove that
$U(k+1)\to U(k)$ is scaled anodyne,
it suffices to show that
$V(k,j+1)\to V(k,j)$ is scaled anodyne
for all $k+1\le j\le n$.

We set $J(k,j)=\{00,\ldots,0k,0j,\ldots,0n,1k,\ldots,1n\}$.
Then we have    
$\tau(k,j)=\Delta^{J(k,j)}$.
The intersection of $\tau(k,j)$
with the underlying simplicial set of $V(k,j+1)$
has the form $\Lambda^{J(k,j)}_{M(k,j)}$,
where
\[ M(k,j)=\left\{\begin{array}{ll}
    \{00,\ldots,0k, 
    1(k+1),\ldots,1n\}&(j=n),\\%[2mm]
    \{00,\ldots,0k, 1(k+1),\ldots,1n,
    0(j+1),\ldots,0n\}&(k+1\le j<n).\\%[2mm]
   \end{array}\right.\]
Hence,
there is a pushout diagram of scaled simplicial sets
\[ \begin{array}{ccc}
     \Lambda^{J(k,j)}_{M(k,j),\S}
%    X(k+1)\cap \Delta^{\{0,\ldots,n,\overline{n},\ldots,\overline{k}\}}
    & \longrightarrow &
      \Delta^{J(k,j)}_{\S}\\
%     \tau(k,k+1)_{\dagger}\\
%    \Delta^{\{0,\ldots,n,\overline{n},\ldots,\overline{k}\}}\\
    \bigg\downarrow & & \bigg\downarrow \\
    V(k,j+1)
    & \longrightarrow & V(k,j).\\
   \end{array}\]
Since 
the top horizontal arrow is scaled anodyne
by Lemma~\ref{lemma:simplex-scaled-anodyne},
so is the bottom horizontal arrow.
\qed

%\subsection{The cartesian fibration
%$\mathbf{p}:\rtwr(Z)\to u_1Z\times u_1Z^{\rm op}$}
\subsection{A model for twisted arrow $\infty$-categories}
\label{subsection:cartesian-fibration-structure}

Let $Z$ be an $\infty$-bicategory.
By Theorem~\ref{thm:twrZ-complete-Segal-space},
$\bbtwr(Z)$
is a complete Segal space.
The inclusion maps
$\Delta^n_{\sharp}\coprod\Delta^{n,{\rm op}}_{\sharp}
\to T(n)$ for $n\ge 0$
induce a map of complete Segal spaces
$q_{\mathbb{T}}: \bbtwr(Z)\to v_1Z\times v_1^{\rm rev}Z$.
By applying the right Quillen equivalence
$i_1^*: \bsetdelta^{\rm CSS}\to \setjoy$,
we obtain a map of 
$\infty$-categories
\[ p_{\mathbb{T}}: \bstwr(Z)\longrightarrow
      u_1Z\times u_1Z^{\rm op},\]
where $\bstwr(Z)$ is an $\infty$-category 
whose set of $n$-simplices
is ${\rm Hom}_{\setsc}(T(n),Z)$.
The goal of this subsection
is to show
%the following proposition.
that $p_{\mathbb{T}}$ is a cartesian fibration
classified by the restricted mapping $\infty$-category
functor,
and hence $\bstwr(Z)$ is a model of
twisted arrow $\infty$-category.
%of the $\infty$-bicategory $Z$.

First,
we will show that
$p_{\mathbb{T}}$ is a categorical fibration
between $\infty$-categories.

\begin{lemma}
The map
$p_{\mathbb{T}}:
\bstwr(Z)\to
      u_1Z\times u_1Z^{\rm op}$
is a categorical fibration
between $\infty$-categories.
\end{lemma}

\proof
First,
we show that 
$q_{\mathbb{T}}:
\bbtwr(Z)\to v_1Z\times v_1^{\rm rev}Z$
is a fibration between fibrant objects
in $\bsetdelta^{\rm CSS}$.
Recall that 
the model category $\bsetdelta^{\rm CSS}$
is a left Bousfield localization
of the Reedy model structure on $\bsetdelta$.
Since $\bbtwr(Z)$ and $v_1Z\times v_1^{\rm rev}Z$
are complete Segal spaces,
in order to show that $q_{\mathbb{T}}$ is a fibration
in $\bsetdelta^{\rm CSS}$,
it suffices to show that
$q_{\mathbb{T}}$ is a Reedy fibration
(cf.~\cite[Proposition~3.3.16]{Hirschhorn}).

For a bisimplicial set $X$,
we denote by $M_nX$ the $n$th matching object.
We set $V= v_1Z\times v_1^{\rm rev}Z$ for simplicity.
We have to show that
the map $f_n: \bbtwr(Z)_n\to M_n\bbtwr(Z)
\times_{M_nV}V_n$ 
is a Kan fibration for all $n\ge 0$.
We observe that
$f_n$ is isomorphic to
the map
$\mapsc(T(n),Z)\to \mapsc(\overline{\partial}T(n),Z)$
induced by
the inclusion map $\overline{\partial}T(n)\to T(n)$,
where $\overline{\partial}T(n)$
is a subcomplex of the scaled simplicial set $T(n)$
given by
\[ \overline{\partial}T(n)=
   \Omega((\Delta^1\times\partial\Delta^n)\cup
   (\partial\Delta^1\times\Delta^n))_{\S}. \]
By Proposition~\ref{prop:cof-induce-map-sc-fibration},
we see that $f_n$ is a Kan fibration. 

Thus,
$q_{\mathbb{T}}$
is a fibration between fibrant objects
in $\bsetdelta^{\rm CSS}$.
By applying the right Quillen equivalence
$i_1^*: \bsetdelta^{\rm CSS}\to \setjoy$,
we see that $p_{\mathbb{T}}$
is a categorical fibration
between $\infty$-categories.
\qed

\bigskip

Next,
we show that $p_{\mathbb{T}}$ is a cartesian fibration.

\begin{proposition}\label{prop:twr-cartesian-fibration}
The map
$p_{\mathbb{T}}: \bstwr(Z)\to u_1Z\times u_1Z^{\rm op}$
is a cartesian fibration.
\end{proposition}

%In the following of this subsection
%we give 
%another proof of Proposition~\ref{prop:twr-cartesian-fibration}
%based on scaled anodyne maps.
In order to prove
Proposition~\ref{prop:twr-cartesian-fibration},
we give a characterization
of $p_{\mathbb{T}}$-cartesian edges.
We define a scaled simplicial set
$\twscncart$ for $n\ge 0$.
The underlying simplicial set of $\twscncart$
is $\Omega^n$
and the set $T_{\twscncart}$ of thin triangles is given by
\[ T_{\twscncart}=
   T_{\twscn}\cup \{\Delta^{\{0i,1(n-1),1n\}}|\ 0\le i<n\}. \]
For a subcomplex $K$ of $\Omega^n$,
we denote by $K_{\rm cart}$
the scaled simplicial set 
whose underlying simplicial set is $K$
equipped with the induced scaling from $\twscncart$.
We denote by $\overline{\Lambda}^n_n\twsccart$
the subcomplex of $\twscncart$
given by
$\Omega(\Delta^1\times\Lambda^n_n
\cup\partial\Delta^1\times\Delta^n)_{\rm cart}$.
We show that
a morphism $\twscone\to Z$
in $\bstwr(Z)$
is a $p_{\mathbb{T}}$-cartesian edge if
it factors through the map
$\twscone\to\twscone_{\rm cart}$.

\begin{lemma}\label{lemma:lambda-n-n-scaled-anodyne}
The map
$\overline{\Lambda}^n_n\twsccart\to
\twscncart$
is scaled anodyne
%for $n\ge 2$.
for $n\ge 1$.
\end{lemma}

\proof
%We consider the case $n\ge 2$.
We recall that
$\sigma(r)=\Delta^{\{00,\ldots,0r,1r,\ldots,1n\}}$
for $0\le r\le n$.
%(see \S\ref{subsection:introduction-twscbullet}).
For $0\le s\le n$,
we define a scaled simplicial set
$X(s)_{\rm cart}$,
which is a subcomplex of
$\Omega^n$
equipped with the induced scaling
from $\twscncart$,
by
\[ X(s)_{\rm cart}=\overline{\Lambda}^n_n\twsccart\cup
        \bigcup_{s\le r\le n} \sigma(r)_{\rm cart}.\]
By definition,
we obtain a filtration
$\overline{\Lambda}^n_n\twsccart = X(n+1)_{\rm cart}
\to\cdots\to X(0)_{\rm cart}=\twscncart$
of scaled simplicial sets
consisting of subcomplexes of $\twscncart$.
In order to prove that
$\overline{\Lambda}^n_n\twsccart\to
\twscncart$ is scaled anodyne,
it suffices to show that
$X(s+1)_{\rm cart}\to X(s)_{\rm cart}$ is scaled anodyne for
all $0\le s\le n$.

We use the notation in the proof of
Lemma~\ref{lemma:inner-twsc}.
We recall that $I(s)=\{00,\ldots,0s,1s,\ldots,1n\}$
for $0\le s\le n$.
The intersection of $\Delta^{I(s)}$
with the underlying simplicial set
of $X(s+1)_{\rm cart}$ has the form
$\Lambda^{I(s)}_{N(s)}$,
where
\[ N(s)= \left\{
   \begin{array}{ll}
    \{1n\}         & (s=0),\\%[2mm]
    \{0s,1n\}    & (0<s<n),\\%[2mm]
    \{0n\}& (s=n).\\
   \end{array}\right.\]
Hence,
there is a pushout diagram of scaled simplicial sets
\[ \xymatrix{
     \Lambda^{I(s)}_{N(s),{\rm cart}}\ar[r]\ar[d]
     & \Delta^{I(s)}_{\rm cart}\ar[d]\\
     X(s+1)_{\rm cart}\ar[r]& X(s)_{\rm cart}.\\
   }\]
By using the fact that
$\Delta^{\{0(n-1),0n,1n\}}$ and $\Delta^{\{0i,1(n-1),1n\}}$
are thin for $0\le i<n$,
we see that the top horizontal arrow is scaled anodyne 
by Lemma~\ref{lemma:simplex-scaled-anodyne}.
Hence,
so is the bottom horizontal arrow.
\qed

\bigskip

We define $\overline{\Lambda}^n_nT'$ by
\[ \overline{\Lambda}^n_nT'=
\overline{\Lambda}^n_nT
\coprod_{\Omega(\Delta^1\times\Delta^{\{n-1,n\}})_{\S}}
  \Omega(\Delta^1\times\Delta^{\{n-1,n\}})_{\rm cart}.\]

\begin{lemma}\label{lemma:lambda-T-prime-scaled-anodyne}
For $n\ge 2$,
the inclusion map
$\overline{\Lambda}^n_nT'
\to \overline{\Lambda}^n_n\twsccart$
is scaled anodyne.
\end{lemma}

\proof
When $n=2$,
it is the identity map.
When $n\ge 3$,
the desired result follows by observing that
$(\Delta^{\{0i,0(n-1),1(n-1),1n\}},T)
\to \Delta^{\{0i,0(n-1),1(n-1),1n\}}_{\sharp}$
is scaled anodyne
for $0\le i< n-1$
by \cite[Remark~3.1.4]{Lurie3},
where $T$ is the set of all $2$-simplices
of $\Delta^{\{0i,0(n-1),1(n-1),1n\}}$
other than $\Delta^{\{0i,1(n-1),1n\}}$.
\qed

\begin{corollary}\label{cor:characterization-p-Cart}
A morphism $\twscone\to Z$ in
$\bstwr(Z)$
is a $p_{\mathbb{T}}$-cartesian edge if 
it factors through the map
$\twscone\to \twscone_{\rm cart}$.
\end{corollary}

\proof
Let $M\in \bstwr(Z)_1$ be a morphism
of the $\infty$-category $\bstwr(Z)$
given by
a map $\overline{M}: \twscone\to Z$ of scaled simplicial sets.
We consider the following commutative diagram
of simplicial sets
\[ \xymatrix{
  \Delta^{\{n-1,n\}}\ar@{^{(}_->}[d]\ar[dr]^{M}& \\
  \Lambda^{n}_n \ar[r]\ar@{^{(}_->}[d] &
  \bstwr(Z)\ar[d]^{p_{\mathbb{T}}} \\
  \Delta^n\ar[r]\ar@{..>}[ur] & u_1Z\times u_1Z^{\rm op} 
}\]
for $n\ge 2$.
By \cite[Remark~2.4.1.4]{Lurie1},
if there exists a lifting
$\Delta^n\to \bstwr(Z)$
that makes the diagram commute,
then $M$ is a $p_{\mathbb{T}}$-cartesian edge.
Unwinding the definitions, 
we see that this is equivalent to the following lifting problem
\[ \xymatrix{
  \Omega(\Delta^1\times\Delta^{\{n-1,n\}})_{\S}
  \ar[dr]^{\overline{M}}\ar@{^{(}_->}[d]& \\
  \overline{\Lambda}^n_nT \ar[r]\ar@{^{(}_->}[d]&
  Z \\
  T(n)\ar@{..>}[ur] & \\
}\]
in $\setsc$.

If $\overline{M}$ factors through $\twscone_{\rm cart}$,
then the map
$\overline{\Lambda}^n_nT\to Z$ extends
to a map $\overline{\Lambda}^n_nT'\to Z$.
We obtain the following commutative diagram
\begin{equation}\label{eq:lambda-n-n-diagram}
  \xymatrix{
  \overline{\Lambda}^n_nT
    \ar[r]\ar[d]&
    \overline{\Lambda}^n_nT'\ar[r]\ar[d]& Z. \\
    T(n) \ar[r]& T(n)_{\rm cart}.\ar@{..>}[ur]& \\    
}\end{equation}
By Lemmas~\ref{lemma:lambda-n-n-scaled-anodyne}
and \ref{lemma:lambda-T-prime-scaled-anodyne},
the middle vertical arrow in (\ref{eq:lambda-n-n-diagram})
is scaled anodyne.
Since $Z$ is a fibrant object in $\setsc$,
the dotted arrow in (\ref{eq:lambda-n-n-diagram})
exists
that makes the right triangle commute.
The desired lifting is obtained as the composite
of $T(n)_{\rm cart}\to Z$ with the inclusion map
$T(n)\to T(n)_{\rm cart}$.
\qed

\proof[Proof of Proposition~\ref{prop:twr-cartesian-fibration}]
We consider a lifting problem depicted
by the following commutative diagram
\begin{equation}\label{eq:p-1-1-lifting-diagram}
  \xymatrix{
    \Delta^{\{1\}}\ar[r]\ar[d]&
    \bstwr(Z)\ar[d]^{p_{\mathbb{T}}}\\
    \Delta^{1}\ar[r]\ar@{..>}[ur]^M&
    u_1Z\times u_1Z^{\rm op}
    \\
} \end{equation}
in $\setdelta$.
%We note that $\Lambda^1_1=\Delta^{\{1\}}\subset \Delta^1$.
Unwinding the definitions,
we see that this is equivalent
to the following lifting problem
\[ \xymatrix{
  \overline{\Lambda}^1_1T \ar[r]\ar[d]& Z\\
  T(1)\ar@{..>}[ur]_{\overline{M}}& \\   
}\]
in $\setsc$.
We would like to have a lifting
$\overline{M}: T(1)\to Z$
which factors through $T(1)_{\rm cart}$.

Since the inclusion map
$\overline{\Lambda}^1_1T\to \overline{\Lambda}^1_1T_{\rm cart}$
is the identity,
we can factor the map
$\overline{\Lambda}^1_1T\to Z$
through $\overline{\Lambda}^1_1T_{\rm cart}$
and we obtain the following commutative diagram
\begin{equation}\label{eq:lambda-1-1-diagram}
  \xymatrix{
   \overline{\Lambda}^1_1T\ar[r]\ar[d]
   & \overline{\Lambda}^1_1T_{\rm cart}\ar[r]\ar[d]
   & Z\\
   T(1)\ar[r]& T(1)_{\rm cart}.\ar@{..>}[ur]& \\
  }\end{equation}
The middle vertical arrow in (\ref{eq:lambda-1-1-diagram})
is scaled anodyne
by Lemma~\ref{lemma:lambda-n-n-scaled-anodyne}.
Since $Z$ is a fibrant object in $\setsc$,
we obtain the dotted arrow $T(1)_{\rm cart}\to Z$
in (\ref{eq:lambda-1-1-diagram})
that makes the right triangle commute.
The composite
$\overline{M}: T(1)\to T(1)_{\rm cart}\to Z$
gives a lifting $M$ in (\ref{eq:p-1-1-lifting-diagram}).
Since $\overline{M}$ factors through $T(1)_{\rm cart}$,
the morphism $M$ is $p_{\mathbb{T}}$-cartesian 
by Corollary~\ref{cor:characterization-p-Cart}.
\qed

\bigskip

%\subsection{Classifying map 
%for the cartesian fibration $\mathbf{p}$}
%\label{subsection:model-twisted-arrow}

%The goal of this subsection is to show the following theorem.
Next,
we show that
$p_{\mathbb{T}}$ is classified
by the restricted mapping $\infty$-category functor.
For this,
it suffices to show that $p_{\mathbb{T}}$
is equivalent to the cartesian fibration
$p: \twr(Z)\to u_1Z\times u_1Z^{\rm op}$.

\begin{theorem}
\label{thm:tw-equivalence-cartesian}
There is an equivalence
of $\infty$-categories
$\twr(Z)\to \bstwr(Z)$
which makes the following diagram commute
\[ \xymatrix{
   \twr(Z)\ar[rr]^{\simeq}\ar[dr]_p&&
   \bstwr(Z)\ar[dl]^{p_{\mathbb{T}}}\\
   & u_1Z\times
     u_1Z^{\rm op}.&\\
}\]  
\end{theorem}

\proof
%By Theorem~\ref{thm:proof-of-equivalnece-model}.
The commutative diagram
is obtained from
the inclusion maps
$\Delta^{\bullet}_{\sharp}\coprod\Delta^{\bullet,{\rm op}}_{\sharp}
\to T(\bullet)\to Q(\bullet)$
of cosimplicial scaled simplicial sets
by applying ${\rm Hom}_{\setsc}(-,Z)$.
The map $T(\bullet)\to Q(\bullet)$
induces an equivalence of complete Segal spaces
$\btwr(Z)\to \bbtwr(Z)$
by Proposition~\ref{prop:proof-of-equivalnece-model}.
By applying the right Quillen equivalence
$i_1^*: \bsetdelta^{\rm CSS}\to \setjoy$,
we see that $\twr(Z)\to \bstwr(Z)$
is an equivalence of $\infty$-categories.
\qed

\begin{corollary}
The cartesian fibration
$p_{\mathbb{T}}: \bstwr(Z)\to u_1Z\times u_1Z^{\rm op}$
is classified by
the restricted mapping $\infty$-category functor
\[ \map{Z}(-,-):
   u_1Z^{\rm op}\times 
   u_1Z\longrightarrow \cat.\]
\end{corollary}

\begin{remark}\rm
There is an isomorphism
$\overline{(-)}: (\Delta^{n}\star\Delta^{n,{\rm op}})^{\rm op}
\to \Delta^{n}\star\Delta^{n,{\rm op}}$
of simplicial sets
given by $\overline{0i}= 1i$
and $\overline{1i}= 0i$.
For a scaled simplicial set $A$ 
whose underlying simplicial set is a subcomplex
of $\Delta^n\star\Delta^{n,{\rm op}}$,
we define its conjugate $A^{\rm conj}$
to be the scaled simplicial set whose underlying
simplicial set is $\overline{A^{\rm op}}$
equipped with a set of thin triangle
$\{\overline{t^{\rm op}}|\ t\in T_A\}$,
where $T_A$ is the set of thin triangles of $A$.
We notice that $Q(n)^{\rm conj}=Q(n)$ for all $n\ge 0$.

In particular,
we have a subcomplex $T(n)^{\rm conj}$ of $Q(n)$.
The underlying simplicial set of
$T(n)^{\rm conj}$ is 
spanned by $(n+1)$-dimensional simplices
\[ \Delta^{\{10,\ldots,1k,0k,\ldots, 0n\}} \]
for $0\le k\le n$.
By assembling $T(n)^{\rm conj}$ for $n\ge 0$,
we obtain a cosimplicial scaled simplicial set 
$T(\bullet)^{\rm conj}=\{T(n)^{\rm conj}\}_{n\ge 0}$
equipped with a map
$T(\bullet)^{\rm conj}\to Q(\bullet)$
of cosimplicial scaled simplicial sets.
For an $\infty$-bicategory $Z$, 
we define
$\bbtwr(Z)^{\rm conj}$
to be the bisimplicial set
$\mapsc(T(\bullet)^{\rm conj},Z)$.

For a monomorphism $A\to B$ of scaled simplicial sets,
where the underlying simplicial sets of $A$ and $B$ 
are subcomplexes of $\Delta^n\star\Delta^{n,{\rm op}}$,
if it is a trivial cofibration of scaled simplicial sets,
then its conjugate $A^{\rm conj}\to B^{\rm conj}$
is also a trivial cofibration. 
By using this fact and 
the functoriality of $(-)^{\rm conj}$,
we see that 
$\bbtwr(Z)^{\rm conj}$
is a complete Segal space
by the conjugates of Lemma~\ref{lemma:twsc-reedy-cofibrant},
Lemma~\ref{lemma:twsc-co-Segal}, 
and Proposition~\ref{prop:proof-of-equivalnece-model}.

We define 
$\bstwr(Z)^{\rm conj}$
to be the $\infty$-category
obtained by applying the right Quillen equivalence
$i_1^*$ to $\bbtwr(Z)^{\rm conj}$.
Although
$(\Delta^{\{00,\ldots,0n\}}_{\sharp})^{\rm conj}
=\Delta^{\{10,\ldots,0n\}}_{\sharp}$
and 
$(\Delta^{\{10,\ldots,0n\}}_{\sharp})^{\rm conj}
=\Delta^{\{00,\ldots,0n\}}_{\sharp}$,
the inclusion map
$\Delta^{\{00,\ldots,0n\}}_{\sharp}
\coprod\Delta^{\{10,\ldots,1n\}}_{\sharp}
\to Q(n)$ is invariant under
$(-)^{\rm conj}$,
and it factors through $T(n)^{\rm conj}$.
The map $\Delta^{\{00,\ldots,0n\}}_{\sharp}
\coprod\Delta^{\{10,\ldots,1n\}}_{\sharp}
\to T(n)^{\rm conj}$ induces a map
$p_{\mathbb{T}}^{\rm conj}:
\bstwr(Z)^{\rm conj}\to u_1Z\times u_1Z^{\rm op}$
of $\infty$-categories.
By the functoriality of $(-)^{\rm conj}$,
we see that
$p_{\mathbb{T}}^{\rm conj}$
is a cartesian fibration
by the conjugates
of Lemmas~\ref{lemma:lambda-n-n-scaled-anodyne}
and \ref{lemma:lambda-T-prime-scaled-anodyne}.

By the conjugate of
Proposition~\ref{prop:proof-of-equivalnece-model},
we have an equivalence
$\btwr(Z)\to \bbtwr(Z)^{\rm conj}$
of complete Segal spaces.
This implies an equivalence
$\twr(Z)\to \bstwr(Z)^{\rm conj}$
of $\infty$-categories.
Since we have monomorphisms
of scaled simplicial sets
$\Delta^{\{00,\ldots,0n\}}_{\sharp}\coprod
\Delta^{\{10,\ldots,1n\}}_{\sharp}\to
T(n)^{\rm conj}\to Q(n)$,
we see that the equivalence 
$\twr(Z)\to \bstwr(Z)^{\rm conj}$
is compatible with the maps $p$
and $p_{\mathbb{T}}^{\rm conj}$
to $u_1Z\times u_1Z^{\rm op}$.
Hence,
we can also regard 
$\bstwr(Z)^{\rm conj}$
as another model of the twisted arrow $\infty$-category.
\end{remark}

%\newpage

\section{Comparisons
with other constructions of twisted arrow categories}
\label{section:comparison}

In this section
we compare the twisted arrow $\infty$-category
in this paper with other related constructions.
In \S\ref{subsection:twisted-2-cell-arrow}
we discuss a relationship with
the construction
of twisted $2$-cell $\infty$-bicategory in \cite{HNP}.
In \S\ref{subsection:other-models-twisted-arrow}
we compare the twisted arrow $\infty$-category
in this paper with the constructions
of twisted arrow $\infty$-category of
an $(\infty,1)$-category
in \cite{Mukherjee-Rasekh}
and \cite{Martini}.

\subsection{Twisted $2$-cell $\infty$-bicategories}
\label{subsection:twisted-2-cell-arrow}

In this subsection
we discuss a relationship between
the construction of twisted arrow
$\infty$-category in this paper 
and that of twisted $2$-cell $\infty$-bicategory of \cite{HNP}.

First,
we recall the notion of
twisted $2$-cell $\infty$-bicategory
introduced in \cite[\S3]{HNP}.
We let $\catdeltaplus$
be the category of $\setplus$-enriched categories
and $\setplus$-enriched functors.
We equip it with the Bergner-Lurie model structure
(cf.~\cite[A.3.2]{Lurie1}).
We refer to an object of $\catdeltaplus$
as a marked simplicial category.
Let $\mathbb{C}$ be a fibrant marked simplicial category.
%that is, a fibrant object of $\catdeltaplus$.
We denote by
${\rm Map}_{\mathbb{C}}:\mathbb{C}^{\rm op}\times\mathbb{C}
\to\setplus$
the mapping marked simplicial set functor.
We note that ${\rm Map}_{\mathbb{C}}$
is a $\setplus$-enriched functor.

By \cite[Theorem~4.2.7]{Lurie3},
there is a Quillen equivalence
\[ \mathfrak{C}^{\rm sc}: \setsc\rightleftarrows
   \catdeltaplus: N^{\rm sc}, \]
where $N^{\rm sc}$
is the scaled nerve functor and
$\mathfrak{C}^{\rm sc}$ is its left adjoint.
We set $\mathbb{D}=\mathbb{C}^{\rm op}\times \mathbb{C}$
and $D=N^{\rm sc}(\mathbb{D})$.
We denote by $\phi: \mathfrak{C}^{\rm sc}(D)
\to \mathbb{D}$
the counit map at $\mathbb{D}$,
which is a weak equivalence in $\catdeltaplus$.

By applying the scaled unstraightening functor
${\rm Un}_{\phi}^{\rm sc}$
(cf.~\cite[\S3.5]{Lurie1} and \cite[\S2.2]{HNP})
to the fibrant $\setplus$-enriched functor
${\rm Map}_{\mathbb{C}}$,
we obtain a scaled cocartesian fibration
in the sense of \cite[Definition~2.7]{HNP}
\[ {\rm Tw}^{\rm sc}(\mathbb{C})\longrightarrow D.\] 
By \cite[Lemma~2.9]{HNP},
the scaled simplicial set
${\rm Tw}^{\rm sc}(\mathbb{C})$ is an $\infty$-bicategory.

We have a twisted arrow construction
on marked simplicial sets
${\rm Tw}^+: \setplus\to \setplus$
(cf.~\cite[the paragraph after Remark~3.2]{HNP}),
which preserves finite products and fibrant objects.
By applying ${\rm Tw}^+$
to the mapping marked simplicial sets of $\mathbb{C}$,
we obtain another fibrant marked simplicial
category $\mathbb{C}_{\rm Tw}$.
The twisted $2$-cell $\infty$-bicategory of $\mathbb{C}$
(\cite[Definition~3.3]{HNP})
is defined by
\[ {\rm Tw}_2(\mathbb{C})={\rm Tw}^{\rm sc}(\mathbb{C}_{\rm Tw}). \]
%We note that 
%thre is a scaled cocartesian fibration
%${\rm Tw}_2(\mathbb{C})\to
%   N^{\rm sc}(\mathbb{C}_{\rm Tw}^{\rm op})\times
%   N^{\rm sc}(\mathbb{C}_{\rm Tw})$.

\if0
We consider $u_1{\rm Tw}_2(\mathbb{C})$
the underlying $\infty$-category of 
the twisted $2$-cell $\infty$-bicategory of $\mathbb{C}$.
Objects of $u_1{\rm Tw}_2(\mathbb{C})$
are $1$-simplices of ${\rm Map}_{\mathbb{C}}(x,y)$
for $x,y\in \mathbb{C}$,
which we think of as morphisms
$(\sigma: f\Rightarrow g)$
between $f,g\in {\rm Map}_{\mathbb{C}}(x,y)$.
Morphisms from
$(\sigma: f\Rightarrow g)$
in ${\rm Map}_{\mathbb{C}}(x,y)$ to
$(\sigma': f'\Rightarrow g')$ in
${\rm Map}_{\mathbb{C}}(x',y')$ are
triples $(\tau_0,\tau_1,\Phi)$,
where $(\tau_0: p_0\Rightarrow q_0)$
is a $1$-simplex in ${\rm Map}_{\mathbb{C}}(x',x)$,
$(\tau_1: p_1\Rightarrow q_1)$ in ${\rm Map}_{\mathbb{C}}(y,y')$,
and $\Phi$ is a $1$-simplex of
${\rm Tw}^+({\rm Map}_{\mathbb{C}}(x',y'))$
such that $\Phi(0)=\tau_1\circ\sigma\circ\tau_0$
and $\Phi(1)=\sigma'$.

Roughly speaking,
a morphism from $\sigma$ to $\sigma'$ is 
and $(\varphi: f'\Rightarrow p_1\circ f\circ p_0)$
and $(\psi: q_1\circ g\circ q_0\Rightarrow g')$
are in ${\rm Map}_{\mathbb{C}}(x',y')$,
such that we have a
$1$-simplex in ${\rm Tw}^+({\rm Map}_{\mathbb{C}}(x',y'))$
depicted as
\[ \xymatrix{
  p_1\circ f\circ p_0
  \ar@{=>}[d]_{\tau_1\circ\sigma\circ\tau_0}
  && f' \ar@{=>}[ll]_-{\varphi}
  \ar@{=>}[d]^{\sigma'}\\
  q_1\circ g\circ q_0\ar@{=>}[rr]_-{\psi}&&
  g'.\\
}\]
\fi

In order to compare ${\rm Tw}_2(\mathbb{C})$
with the construction of twisted arrow
$\infty$-category in this paper,
we will construct a functor
${\rm Tw}_2(\mathbb{C})\to {\rm Tw}^{\rm sc}(\mathbb{C})$
of $\infty$-bicategories.
For any fibrant marked simplicial set $W$,
we have a marked left fibration
${\rm Tw}^+(W)\to W^{\rm op}\times W$
in the sense of \cite[Definition~2.16]{HNP}.
By composing with the projection
$W^{\rm op}\times W\to W$,
we obtain a map of marked simplicial sets
${\rm Tw}^{+}(W)\to W$.
By applying this construction
to the mapping scaled simplicial sets of $\mathbb{C}$,
we obtain 
a $\setplus$-enriched functor
$\mathbb{C}_{\rm Tw}\to\mathbb{C}$.
We set $\mathbb{D}_{\rm Tw}=\mathbb{C}_{\rm Tw}^{\rm op}\times
\mathbb{C}_{\rm Tw}$ and
$D_{\rm Tw}=N^{\rm sc}(\mathbb{D}_{\rm Tw})$.
The functor $\mathbb{C}_{\rm Tw}\to\mathbb{C}$
induces an $\setplus$-enriched
functor $\psi: \mathbb{D}_{\rm Tw}\to\mathbb{D}$,
a map of scaled simplicial sets
$f=\mathfrak{C}^{\rm sc}(\psi): D_{\rm Tw}\to D$, 
and a natural transformation
of $\setplus$-enriched functors
${\rm Map}_{\mathbb{C}_{\rm Tw}}\to
{\rm Map}_{\mathbb{C}}\circ \psi$.
By \cite[Remarks~3.5.16 and 3.5.17]{Lurie3}, 
we obtain the following commutative diagram
\[ \xymatrix{
     {\rm Tw}_2(\mathbb{C})\ar[r]\ar[d]&
%     f^*{\rm Tw}^{\rm sc}(\mathbb{C})\ar[r]\ar[d]&
     {\rm Tw}^{\rm sc}(\mathbb{C})\ar[d]  \\
     D_{\rm Tw}\ar[r]^f&
%    D_{\rm Tw}\ar[r]_f&
     D. \\
}\]
%where the right square is a pullback diagram.

Although the twisted $2$-cell $\infty$-bicategory
${\rm Tw}_2(\mathbb{C})$ is more general
than ${\rm Tw}^{\rm sc}(\mathbb{C})$,
the construction of ${\rm Tw}^{\rm sc}(\mathbb{C})$
is more directly connected to
that of twisted arrow $\infty$-category
in this paper.
We will compare the construction
of twisted arrow $\infty$-category
for an $\infty$-bicategory in this paper
and that of ${\rm Tw}^{\rm sc}(-)$.

Let $Z$ be an $\infty$-bicategory.
We take a fibrant replacement
$\mathfrak{C}^{\rm sc}(Z)\to \mathbb{C}(Z)$
in $\catdeltaplus$.
%We set $\mathcal{D}(Z)=\mathbb{C}(Z)^{\rm op}\times
%\mathbb{C}(Z)$.
%and
%$D(Z)=N^{\rm sc}(\mathcal{D}(Z))$.
%The fibrant replacement
%$\mathfrak{C}^{\rm sc}(Z)\to \mathbb{C}(Z)$
It induces a fibrant replacement
$\varphi: \mathfrak{C}^{\rm sc}(Z^{\rm op}\times Z)
\to \mathbb{C}(Z)^{\rm op}\times\mathbb{C}(Z)$.
By applying the unstraightening functor
${\rm Un}_{\varphi}^{\rm sc}$
to the fibrant $\setplus$-enriched functor
${\rm Map}_{\mathbb{C}(Z)}$,
we obtain a scaled cocartesian fibration
of $\infty$-bicategories
${\rm Tw}^{\rm sc}(\mathbb{C}(Z))\to
   Z^{\rm op}\times Z$.
Furthermore,
by applying the right Quillen functor
$u_1: \setsc\to \setjoy$,
we obtain a cocartesian fibration
of $\infty$-categories
\[ {\rm Tw}^l(\mathbb{C}(Z))
   \longrightarrow u_1Z^{\rm op}\times u_1Z,\]
where ${\rm Tw}^l(\mathbb{C}(Z))=
u_1{\rm Tw}^{\rm sc}(\mathbb{C}(Z))$
is the underlying $\infty$-category
of the $\infty$-bicategory
${\rm Tw}^{\rm sc}(\mathbb{C}(Z))$.
By construction,
it is classified by
the restricted mapping $\infty$-category functor
${\rm Map}_Z: u_1Z^{\rm op}\times u_1 Z
   \to\cat$.

We say that a cartesian fibration
and a cocartesian fibration
of $\infty$-categories are dual to each other
if they classify the same functor to $\cat$.
A model of dual (co)cartesian fibrations 
was constructed in \cite{BGN}.
Since the cartesian fibration
$\bstwr(Z)\to u_1Z^{\rm op}\times u_1Z$
of the twisted arrow $\infty$-category
for an $\infty$-bicategory $Z$
in this paper is classified
by the restricted mapping $\infty$-category
functor
${\rm Map}_Z: u_1Z^{\rm op}\times u_1 Z\to\cat$,
it is a dual to the cocartesian fibration
${\rm Tw}^l(\mathbb{C}(Z))\to u_1Z^{\rm op}\times u_1Z$.

\if0
it is desirable to have
a construction of twisted arrow category

Actually,
it is a restriction
of a notion of twisted $2$-cell $\infty$-bicategory.
We start with a fibrant $\setplus$-enriched category $\mathbb{C}$,
which can be regarded a model of an $\infty$-bicategory.
\fi

%\newpage

\subsection{Other models of twisted arrow $\infty$-categories}
\label{subsection:other-models-twisted-arrow}

In this subsection
we compare other constructions
of twisted arrow $\infty$-categories.
%given in \cite{Mukherjee-Rasekh}
%and \cite{Martini}.

\if0
First,
we recall 
there are two Quillen equivalences between $\bsetdelta^{\rm CSS}$
and $\setjoy$
by Joyal-Tierney~\cite{Joyal-Tierney},
where $\bsetdelta^{\rm CSS}$ is 
the category of bisimplicial sets 
equipped with complete Segal space model structure,
and $\setjoy$ is the category of simplicial sets
with Joyal model structure.
Let $p_1: \Delta^{\rm op}\times\Delta^{\rm op}\to \Delta^{\rm op}$
be the first projection,
and let $i_1: \Delta^{\rm op}\simeq \{[0]\}\times\Delta^{\rm op}
\to \Delta^{\rm op}\times\Delta^{\rm op}$ be the inclusion. 
These maps induce an adjunction
\[ p_{1}^*: \setjoy\rightleftarrows \bsetdelta^{\rm CSS}: i_1^*,\]
which is a Quillen equivalence by
\cite[Theorem~4.11]{Joyal-Tierney}.

We recall that $E(n)$
is the nerve of the groupoid
$[n]^{\rm gpd}$
freely generated by the category $[n]$.
Let $t: \Delta\times\Delta\to \setdelta$
be the functor given by
$t([m],[n])=\Delta^m\times E(n)$.
The functor $t_!: \bsetdelta\to \setdelta$
is a left Kan extension of $t$
along the Yoneda embedding $\Delta\times\Delta\to \bsetdelta$.
The functor $t_!$ admits a right adjoint
$t^!: \bsetdelta\to \setdelta$ given by
\[ t^!(X)_{m,n}={\rm Hom}_{\setdelta}
   (\Delta^m\times E(n), X). \]
The adjunction 
\[ t_!: \bsetdelta^{\rm CSS}\rightleftarrows \setjoy: t^!. \]
is a Quillen equivalence by \cite[Theorem~4.12]{Joyal-Tierney}.
\fi

We consider relationships between
various constructions of twisted arrow $\infty$-categories.
First,
we compare 
the construction by Abell\'{a}n Garc\'{i}a and Stern
with that of Lurie.
Abell\'{a}n Garc\'{i}a and Stern generalized
Lurie's construction of twisted arrow $\infty$-category
of an $(\infty,1)$-category to that of an $(\infty,2)$-category.
The functor
$(-)_{\sharp}: \setjoy\to \setsc$
is also right Quillen
(cf.~\cite[(2.3)]{HNP}),
which induces a functor
$(-)_{\sharp}: (\setjoy)^{\circ}\to (\setsc)^{\circ}$.
%$\cat\to \bicat$
%of $\infty$-categories.
We have the following commutative diagram
\[ \xymatrix{
  (\setjoy)^{\circ}\ar[rrr]^{{\rm Tw}^{\rm Lurie}}\ar@{=}[d]&&&
  (\setjoy)^{\circ}\ar@{=}[d]\\
  (\setjoy)^{\circ}\ar[r]^{(-)_{\sharp}}&
  (\setsc)^{\circ}\ar[rr]^{{\rm Tw}^{\rm AGS}}&&
  (\setjoy)^{\circ},\\
}\]
where ${\rm Tw}^{\rm Lurie}$
is Lurie's twisted arrow construction,
and ${\rm Tw}^{\rm AGS}$
is Abell\'{a}n Garc\'{i}a and Stern's.

Next,
we compare the functor
$\btwr: (\setsc)^{\circ}\to (\bsetdelta^{\rm CSS})^{\circ}$
constructed in \S\ref{section:twisted-arrow-categories}
with the functor
${\rm Tw}^{\rm AGS}: (\setsc)^{\circ}\to
  (\setjoy)^{\circ}$.
%by Abell\'{a}n Garc\'{i}a and Stern.
The functor
$\btwr$
is a lifting of ${\rm Tw}^{\rm AGS}$
in the sense that we have the following commutative diagram
\[ \xymatrix{
  (\setsc)^{\circ}\ar[rr]^{\btwr}&&
  (\bsetdelta^{\rm CSS})^{\circ}\ar[d]^{i_1^*}\\
  (\setsc)^{\circ}\ar[rr]^{{\rm Tw}^{\rm AGS}}\ar@{=}[u]
  && (\setjoy)^{\circ}.\\
}\]

Now,
we compare
the construction of twisted arrow $\infty$-category
by Mukherjee and Rasekh \cite{Mukherjee-Rasekh}
with $\btwr$.
Let $\epsilon: \Delta\to\Delta$
be a functor given by $[n]\mapsto [n]\star[n]^{\rm op}$.
We consider a functor
${\rm Tw}^{\rm MR}: \bsetdelta\to\bsetdelta$
given by
${\rm Tw}^{\rm MR}(W)=W\circ\epsilon^{\rm op}$
for $W\in \bsetdelta$,
where we regard $\bsetdelta$ as
${\rm Fun}(\Delta^{\rm op},\setdelta)$.
This is a variant of the construction
by Mukherjee and Rasekh \cite{Mukherjee-Rasekh},
which is given by $W\mapsto W\circ \epsilon'$
with $\epsilon'([n])=[n]^{\rm op}\star[n]$.
By using \cite[Theorem~3.8]{Mukherjee-Rasekh},
we see that ${\rm Tw}^{\rm MR}$ determines a functor
${\rm Tw}^{\rm MR}: (\bsetdelta^{\rm CSS})^{\circ}
 \to(\bsetdelta^{\rm CSS})^{\circ}$.
%by \cite[Definition~4.2.4 and Proposition~4.2.5]{Martini}.
We can verify that 
we have the following commutative diagram
\[ \xymatrix{
  (\setjoy)^{\circ}\ar[rrr]^{(-)_{\sharp}}
  \ar[d]_{t^!}&&&
  (\setsc)^{\circ}\ar[d]^{\btwr}\\
  (\bsetdelta^{\rm CSS})^{\circ}
  \ar[rrr]^{{\rm Tw}^{\rm MR}}\ar[d]_{i_1^*}
  &&& (\bsetdelta^{\rm CSS})^{\circ}
  \ar[d]^{i_1^*}\\
  (\setjoy)^{\circ}\ar[r]^{(-)_{\sharp}}&
  (\setsc)^{\circ}\ar[rr]^{\mathrm{Tw}^{\rm AGS}}&&
  (\setjoy)^{\circ}.\\
}\]
%and 
%\[ \xymatrix{
%  (\setjoy)^{\circ}\ar[r]^{(-)_{\sharp}}&
%  (\setsc)^{\circ}\ar[rr]^{\mathrm{Tw}^{\rm AGS}}&&
%  (\setjoy)^{\circ}\\
%  (\bsetdelta^{\rm CSS})^{\circ}\ar[rrr]^{\mathrm{Tw}^{\rm MR}}
%  \ar[u]^{i_1^*}_{\simeq}&&&
%  (\bsetdelta^{\rm CSS})^{\circ}\ar[u]_{i_1^*}^{\simeq}.\\
%}\]

Finally,
Martini \cite{Martini}
has developed the theory of categories internal
to an $\infty$-topos $\mathcal{B}$. 
He has introduced
a notion of $\mathcal{B}$-category.
This is a generalization of complete Segal space
since the notion of $\mathcal{B}$-category
coincides with that of complete Segal space
when $\mathcal{B}$ is the $\infty$-category
of $\infty$-groupoids $\mathcal{S}$.
We denote by ${\rm Cat}(\mathcal{B})$
the $\infty$-category 
of $\mathcal{B}$-categories.
Thus,
there is an equivalence of $\infty$-categories
between ${\rm Cat}(\mathcal{S})$ and 
the underlying $\infty$-category
%$(\bsetdelta^{\rm CSS})_{\infty}$
of $\bsetdelta^{\rm CSS}$.
He has also introduced a construction
of twisted arrow $\infty$-category for
a $\mathcal{B}$-category,
which determines a functor
${\rm Tw}^{\rm Mar}_{\mathcal{B}}: {\rm Cat}(\mathcal{B})\to
{\rm Cat}(\mathcal{B})$
by \cite[Definition~4.2.4 and Proposition~4.2.5]{Martini}.

We compare the functor ${\rm Tw}^{\rm Mar}_{\mathcal{S}}$
with ${\rm Tw}^{\rm MR}$ and $\btwr$.
For this purpose,
we introduce some notation.
For a model category $\mathcal{M}$,
we denote by $\mathcal{M}_{\infty}$
the underlying $\infty$-category.
Suppose that $F: \mathcal{M}^{\circ}\to \mathcal{N}^{\circ}$
is a functor between full subcategories
of model categories $\mathcal{M}$ and $\mathcal{N}$
spanned by fibrant-cofibrant objects
that preserves weak equivalences.
We denote by $F_{\infty}: \mathcal{M}_{\infty}\to\mathcal{N}_{\infty}$
the induced functor between the underlying $\infty$-categories.

By the construction
of ${\rm Tw}^{\rm Mar}_{\mathcal{S}}$
in \cite[Definition~4.2.4]{Martini},
%${\rm Tw}^{\rm Mar}_{\mathcal{S}}$
it is equivalent to the functor ${\rm Tw}^{\rm MR}_{\infty}$
of $\infty$-categories
induced by ${\rm Tw}^{\rm MR}$.
Therefore,
we have the following commutative diagram of $\infty$-categories
\[ \xymatrix{
  (\setjoy)_{\infty}\ar[rr]^{(-)_{\sharp,\infty}}
  \ar[d]_{t^!_{\infty}}&&
  (\setsc)_{\infty}\ar[d]^{\btwr_{\infty}}\\
  (\bsetdelta^{\rm CSS})_{\infty}\ar@{=}[d]
  \ar[rr]^{\mathrm{Tw}^{\rm MR}_{\infty}}
  && (\bsetdelta^{\rm CSS})_{\infty}\ar@{=}[d]\\
  {\rm Cat}(\mathcal{S})\ar[rr]^{{\rm Tw}^{\rm Mar}_{\mathcal{S}}}&&
  {\rm Cat}(\mathcal{S}).\\
}\]

\if0
Let $\Gamma: \mathcal{B}\to \mathcal{S}$
be a unique geometric morphism,
which induces a functor
$\Gamma: {\rm Cat}(\mathcal{B})\to{\rm Cat}(\mathcal{S})$.
We have the following commutative diagram
\[ \xymatrix{
   {\rm Cat}(\mathcal{B})\ar[rr]^{{\rm TW}^{\rm Mar}}
   \ar[d]_{\Gamma}&& {\rm Cat}(\mathcal{B})\ar[d]^{\Gamma}\\
   {\rm Cat}(\mathcal{S})\ar[rr]^{{\rm Tw}^{\rm MR}}&&
   {\rm Cat}(\mathcal{S}).\\
}\]
\fi

%\newpage

\section{An extension of ${\rm TW}^r(-)$ to a right Quillen functor}
\label{section:extension-TWr}

In this section
we construct a Quillen adjunction
\[ L_{Q,E}: \bsetdelta^{\rm CSS}\rightleftarrows \setsc: R_{Q,E} \]
whose right adjoint is equivalent to 
the twisted arrow construction
${\rm TW}^r(-): (\setsc)^{\circ}\to (\bsetdelta^{\rm CSS})^{\circ}$
when restricted to the full subcategory $(\setsc)^{\circ}$
of $\infty$-bicategories.

First,
we construct the functor $L_{Q,E}: \bsetdelta\to\setsc$.
We recall that $Q(\bullet)$ is a cosimplicial 
scaled simplicial set
which is Reedy cofibrant by Lemma~\ref{lemma:twr-Z-Segal-space}.
We also recall that
the simplicial set $E(n)$ is
the nerve of the groupoid
$[n]^{\rm gpd}$
freely generated by the category $[n]$.
The assignment $[n]\mapsto E(n)$
determines a cosimplicial simplicial set $E(\star)$.
We will show that $E(\star)_{\sharp}$ is a Reedy cofibrant
cosimplicial scaled simplicial set.

\begin{lemma}
The cosimplicial simplicial set
$E(\star)$ is Reedy cofibrant.
\end{lemma}

\proof
This follows by observing that 
the $n$th latching object $L_nE(\star)$
is isomorphic to the subcomplex 
$\cup_{i=0}^n N(([n]-\{i\})^{\rm gpd})$
of $E(n)$.
\qed

\bigskip

By \cite[Proposition~3.1.5]{Hovey},
the cosimplicial object $E(\star)$
induces an adjunction
\[ L_E: \setkan\rightleftarrows \setjoy : R_E, \]
where $L_E(K)$ is given by
the left Kan extension of $E(\star): \Delta\to \setdelta$
along the Yoneda embedding $\Delta\to \setdelta$.
%the coend
%$E(\star)\otimes K$ for $K\in \setdelta$.
For notational ease,
we write $E(K)$
for the simplicial set $L_E(K)$.
By \cite[Theorem~1.19]{Joyal-Tierney}
(\cite[Theorem~6.22]{Joyal-1}),
the adjoint pair $(L_E,R_E)$ is a Quillen adjunction.
%By \cite[Proposition~5.4.1]{Hovey},
%the left adjoint $L_E: \setdelta\to\setsc$
%preserves cofibrations.

Since the functor $(-)_{\sharp}: \setjoy\to \setsc$
is a left Quillen functor,
we obtain the following corollary.

\begin{corollary}
The cosimplicial scaled simplicial set
$E(\star)_{\sharp}$ is Reedy cofibrant.
\end{corollary}

By the cartesian product of the cosimplicial scaled simplicial sets
$Q(\bullet)$ and $E(\star)_{\sharp}$, 
we obtain a functor $Q(\bullet)\times E(\star)_{\sharp}:
\Delta\times\Delta\to \setsc$ 
which assigns to $([m],[n])$
the scaled simplicial set $Q(m)\times E(n)_{\sharp}$.
We define a functor $L_{Q,E}: \bsetdelta\to \setsc$
to be the left Kan extension of $Q(\bullet)\times E(\star)_{\sharp}$
along the Yoneda embedding
$\Delta\times\Delta\to \bsetdelta$.

%\begin{definition}\rm
%The functor $L_{Q,E}$ admits a right adjoint
%and we denote it by $R_{Q,E}$.
%\end{definition}

The functor $L_{Q,E}$ has a right adjoint
$R_{Q,E}: \setsc\to \bsetdelta$.
Let $\pi_i: \Delta^{\rm op}\times
\Delta^{\rm op}\to\Delta^{\rm op}$
be the $i$th projection for $i=1,2$. 
We define an object $F(m,n)$ of $\bsetdelta$ 
to be $\pi_1^*(\Delta^m)\times \pi_2^*(\Delta^n)$.
We note that there is a natural isomorphism
${\rm Hom}_{\bsetdelta}(F(m,n),A)\cong A_{m,n}$
for any $A\in \bsetdelta$.
%The functor $L_{Q,E}$ has a right adjoint
%$R_{Q,E}: \setsc\to \bsetdelta$
%\begin{definition}\rm
For any scaled simplicial set $X$,
the bisimplicial set $R_{Q,E}(X)$ satisfies
the following natural bijection
\[ {\rm Hom}_{\bsetdelta}(F(m,n),R_{Q,E}(X))
   \cong
{\rm Hom}_{\setsc}(Q(m)\times E(n)_{\sharp},X )\]
for any $m,n\ge 0$.
%\end{definition}

Next,
we will show that $R_{Q,E}(Z)$ is a complete Segal space
for any $\infty$-bicategory $Z$.  

\begin{lemma}
For any $\infty$-bicategory $Z$,
the bisimplicial set $R_{Q,E}(Z)$ is Reedy fibrant. 
\end{lemma}

\proof
By Lemma~\ref{lemma:twr-Z-Segal-space},
${\rm FUN}(Q(\bullet), Z)$ is Reedy fibrant
in the category of simplicial object of $\setsc$.
By applying the right Quillen functor $u_1$,
${\rm Fun}(Q(\bullet), Z)$ is Reedy fibrant
in the category of simplicial object of $\setjoy$.
We notice that there is a natural isomorphism
of bisimplicial sets
$R_E{\rm Fun}(Q(\bullet,Z))\cong
R_{Q,E}(Z)$.
The lemma follows from
the fact that
$R_E: \setjoy\to \setdelta^{\rm Kan}$
is a right Quillen functor
by \cite[Theorem~1.19]{Joyal-Tierney}
(\cite[Theorem~6.22]{Joyal-1}).
\qed

\if0
\begin{lemma}
For any $\infty$-category $Z$,
$R_{Q,E}(Z)$ is a Segal space.  
\end{lemma}

\proof
By Proof of Proposition~\ref{prop:lifting-SS},
the cosimplicial scaled simplicial set
$Q(\bullet)$ satisfies the co-Segal condition.
This implies that
${\rm Fun}(Q(\bullet),Z)$ satisfies
the Segal condition in the model category
$\setjoy$.
The lemma follows by applying
the right Quillen functor
$R_E: \setjoy\to \setdelta^{\rm Kan}$
to ${\rm Fun}(Q(\bullet),Z)$.
\qed
\fi

\begin{lemma}\label{lemma:RQE-TWr-levelwise-eq}
For any $\infty$-bicategory $Z$,
there is a natural map of bisimplicial sets
$R_{Q,E}(Z)\to  {\rm TW}^r(Z)$,
which is a levelwise trivial Kan fibration.
%in $\bsetdelta^{\rm Reedy}$,
%that is,
%a levelwise trivial Kan fibration.
\end{lemma}

\proof
For a scaled simplicial set $Z$,
there is a natural isomorphism of bisimplicial sets
\[ {\rm Hom}_{\setsc}(Q(\bullet)\times E(\star)_{\sharp}, Z)
   \cong
   {\rm Hom}_{\setdelta}(E(\star), {\rm Fun}(Q(\bullet),Z)).\]
When $Z$ is an $\infty$-bicategory,
we have a natural isomorphism of bisimplicial sets
\[ %\begin{array}{rcl}
   {\rm Hom}_{\setdelta}
   (\Delta^{\star},\mapsc(Q(\bullet),Z))\cong
%   {\rm Hom}_{\setdelta}(\Delta^{\star},
%   u_0{\rm Fun}(Q(\bullet),X)^{\natural})\\
%   &\cong&
   {\rm Hom}_{\setplus}((\Delta^{\star})^{\sharp},
   {\rm Fun}(Q(\bullet),Z)^{\natural}).
   %\end{array}
\]
The proposition follows from the fact that
the inclusion map $\Delta^{\star}\to E(\star)$
induces a natural trivial Kan fibration
${\rm Hom}_{\setdelta}(E(\star),Y)
   \to
   {\rm Hom}_{\setdelta}(\Delta^{\star},Y^{\simeq}) 
   \cong
   {\rm Hom}_{\setplus}((\Delta^{\star})^{\sharp},Y^{\natural})$
for any $\infty$-category $Y$
by \cite[Proposition~1.20]{Joyal-Tierney}
(\cite[Proposition~6.26]{Joyal-1}).       
\if0
we have the following isomorphisms
\[ \begin{array}{rcl}
    R_{Q,E}(X)_{m,n}&\cong&
    {\rm Hom}_{\setsc}(Q(m)\times E(n)_{\sharp},X)\\
    &\cong&
    {\rm Hom}_{\setsc}(E(n)_{\sharp},{\rm FUN}(Q(m),X))\\
    &\cong&
    {\rm Hom}_{\setdelta}(E(n), {\rm Fun}(Q(m),X))\\
    &\cong&
    {\rm Hom}_{\setplus}((\Delta^n)^{\sharp},{\rm Fun}(Q(m),X)^{\natural})\\
    &\cong&
    {\rm Hom}_{\setdelta}(\Delta^n,u_0{\rm Fun}(Q(m),X)^{\natural})\\
    &\cong&
    {\rm Map}^{\rm sc}(Q(m),X)_n.\\ 
   \end{array}\]
This implies that there is a natural isomorphism
$R_{Q,E}(X)\cong {\rm TW}^r(X)$.
\fi
\qed

\begin{corollary}\label{cor:RQE-complete-Segal}
For any $\infty$-bicategory $Z$,
$R_{Q,E}(Z)$ is a complete Segal space. 
\end{corollary}

\proof
Since ${\rm TW}^r(Z)$ is a complete Segal space
by Proposition~\ref{prop:lifting-SS} and
Theorem~\ref{thm:GHmodel-complete-Segal-TwZ},
the corollary follows from
Lemma~\ref{lemma:RQE-TWr-levelwise-eq}.
\qed

\bigskip

By Lemma~\ref{lemma:RQE-TWr-levelwise-eq},
the functor $R_{Q,E}$ is equivalent to ${\rm TW}^r(-)$
when restricted to $(\bsetdelta^{\rm CSS})^{\circ}$.
Next,
we will show that 
$R_{Q,E}$ is a right Quillen functor.

\begin{theorem}\label{thm:extension-TW-right-quillen}
The adjunction $(L_{Q,E},R_{Q,E})$ induces a Quillen adjunction
\[ L_{Q,E}: \bsetdelta^{\rm CSS}\rightleftarrows \setsc: R_{Q,E}. \] 
\end{theorem}

In order to prove
Theorem~\ref{thm:extension-TW-right-quillen},
first, we will show that
the left adjoint $L_{Q,E}$ preserves cofibrations.
The complete Segal space model structure
on the category of bisimplicial sets
is a left Bousfield localization
of the Reedy model structure.
We denote by $\bsetdelta^{\rm Reedy}$
the category of bisimplicial sets
equipped with Reedy model structure.
The Reedy model structure on $\bsetdelta$
is cofibrantly generated.
The pushout product of
$\pi_1^*(\partial\Delta^m)\to \pi_1^*(\Delta^m)$
and
$\pi_2^*(\partial\Delta^m)\to \pi_2^*(\Delta^m)$,
we obtain a map of bisimplicial sets
$\partial F(m,n)\to F(m,n)$,
where 
$\partial F(m,n)=
(\pi_1^*(\partial\Delta^m)\times \pi_2^*(\Delta^n))
\coprod_{\pi_1^*(\partial \Delta^m)\times \pi_2^*(\partial\Delta^n)}
(\pi_1^*(\Delta^m)\times \pi_2^*(\partial \Delta^n))$.
We can take a generating set of cofibrations
of $\bsetdelta^{\rm Reedy}$ as
\[ \{\partial F(m,n)\to F(m,n)|\ m,n\ge 0\}\]
(cf.~\cite[Theorem~15.6.27]{Hirschhorn}).

\if0
\begin{lemma}
Let $f: X\to Y$ be a fibration between $\infty$-bicategories
in $\setsc$.
The map $R_{Q,E}(f): R_{Q,E}(X)\to R_{Q,E}(Y)$
is a fibration between complete Segal spaces
in $\bsetdelta^{\rm CSS}$.
\end{lemma}

\proof
\qed
\fi

\begin{lemma}\label{lemma:LQE-preserve-cofibrations}
The left adjoint $L_{Q,E}$ preserves cofibrations.
\end{lemma}

\proof
It suffices to show that
$L_{Q,E}\partial F(m,n)\to L_{Q,E}F(m,n)$
is a monomorphism for any $m,n\ge 0$.
We have an isomorphism of scaled simplicial
sets $L_{Q,E}F(m,n)\cong Q(m)\times E(n)_{\sharp}$.
The lemma follows by observing that
the map $L_{Q,E}\partial F(m,n)\to F(m,n)$ is isomorphic to
the pushout product of monomorphisms
$Q(\partial\Delta^m)\to Q(m)$ and
$E(\partial \Delta^n)_{\sharp}\to E(n)_{\sharp}$,
which is a monomorphism.
%$\partial (Q(m)\times E(n)_{\sharp})\to (Q(m)\times E(n)_{\sharp})$,
%where 
%$\partial (Q(m)\times E(n)_{\sharp})=
%Q(Q(\partial \Delta^m)\times E(\Delta^n)_{\sharp})
%\coprod_{Q(\partial \Delta^m)\times E(\partial \Delta^n)_{\sharp}}
%(Q(\Delta^m)\times E(\partial \Delta^n))_{\sharp}$
%of $Q(m)\times E^n_{\sharp}$.
\qed

\bigskip

By \cite[Proposition~7.15]{Joyal-Tierney}
(\cite[Proposition~E.2.14]{Joyal-1})
and Lemma~\ref{lemma:LQE-preserve-cofibrations},
in order to prove that $(L_{Q,E},R_{Q,E})$
is a Quillen adjunction,
it suffices to show that
the right adjoint $R_{Q,E}$ preserves
fibrations between fibrant objects.

\begin{lemma}\label{lemma:RQE-preserve-fibrations}
The right adjoint $R_{Q,E}$ preserves
fibrations between fibrant objects.
\end{lemma}

\proof
We write $R(Z)$ for $R_{Q,E}(Z)$ for simplicity.
Let $f: X\to Y$ be a fibration of scaled simplicial sets
between $\infty$-bicategories.
By Corollary~\ref{cor:RQE-complete-Segal},
$R(X)$ and $R(Y)$
are complete Segal spaces and hence
they are fibrant objects in $\bsetdelta^{\rm CSS}$.
Since the model category $\bsetdelta^{\rm CSS}$
is a left Bousfield localization
of $\bsetdelta^{\rm Reedy}$,
it suffices to show that the map
$R(f): R(X)\to R(Y)$ is a Reedy fibration
by \cite[Proposition~3.3.16(1)]{Hirschhorn}.

%$F(k)_{m,n}={\rm Hom}_{\Delta}([m],[k])$.
The model category $\bsetdelta^{\rm Reedy}$
is a simplicial model category
(cf.~\cite[Theorem~15.3.4(3)]{Hirschhorn}).
We denote by
${\rm Map}^{\rm bs}(A,B)$
the mapping simplicial set of $A,B\in \bsetdelta$. 
We set $F(k)=F(k,0)$ and $F(\partial \Delta^k)=\partial F(k,0)$
for simplicity.
The inclusion map $i: F(\partial\Delta^k)\to F(k)$
induces a map
\[ (i,f)^{\rm bs}: {\rm Map}^{\rm bs}(F(k), R(X))\to
   {\rm Map}^{\rm bs}(F(\partial \Delta^k),R(X))
   \times_{{\rm Map}^{\rm bs}(F(\partial\Delta^k),R(Y))}
         {\rm Map}^{\rm bs}(F(k),R(Y)) \]
of simplicial sets.
In order to show that the map
$R(f)$ is a Reedy fibration,
it suffices to show that 
the map $(i,f)^{\rm bs}$ is a Kan fibration
by \cite[Proposition~9.4.4(3) and
Lemma~9.4.7]{Hirschhorn}.

Since $\setsc$ is a cartesian closed model category
and ${\rm FUN}(-,-)$ is the mapping object,
we have a fibration of $\infty$-bicategories
\[ {\rm FUN}(i,f): {\rm FUN}(Q(k),X)\to
   {\rm FUN}(Q(\partial\Delta^k),X)\times_{{\rm FUN}(Q(\partial\Delta^k),Y)}
   {\rm FUN}(Q(k),Y).\]
By applying the right Quillen functor
$u_1: \setsc\to\setjoy$,
we obtain a categorical fibration
of $\infty$-categories
\[ {\rm Fun}(i,f): {\rm Fun}(Q(k),X)\to
   {\rm Fun}(Q(\partial\Delta^k),X)\times_{{\rm Fun}(Q(\partial\Delta^k),Y)}
   {\rm Fun}(Q(k),Y).\]
We notice that there are natural isomorphisms
of Kan complexes
%\[ R_E{\rm Fun}(Q(k),Z)\cong
%   {\rm Map}^{\rm bs}(F(k),R(Z)),
%   \quad
%   R_E{\rm Fun}(Q(\partial\Delta^k),R(Z))
%   \cong
%   {\rm Map}^{\rm bs}(F(\partial \Delta^k),Z)
%\]
%
%\bigskip
%  
\[ \begin{array}{rcl}
    R_E{\rm Fun}(Q(k),Z)&\cong&
   {\rm Map}^{\rm bs}(F(k),R(Z)),\\
   \quad
   R_E{\rm Fun}(Q(\partial\Delta^k),Z)
   &\cong&
   {\rm Map}^{\rm bs}(F(\partial \Delta^k),R(Z))\\
   \end{array}\]
for any $\infty$-bicategory $Z$.
By applying the right Quillen functor
$R_E: \setjoy\to \setkan$ to ${\rm Fun}(i,f)$,
we see that $(i,f)^{\rm bs}$
is a Kan fibration 
between Kan complexes.
\qed

\proof[Proof of Theorem~\ref{thm:extension-TW-right-quillen}]
The theorem follows from
Lemmas~\ref{lemma:LQE-preserve-cofibrations}
and \ref{lemma:RQE-preserve-fibrations}
by using \cite[Proposition~7.15]{Joyal-Tierney}
(\cite[Proposition~E.2.14]{Joyal-1}).
\qed

%\newpage

%\input{ref}

%}

\end{document}